%% file: ms.tex
\documentclass[a4paper]{amsart}

\input{head}
\usepackage[toc,page]{appendix}
\usepackage{etoolbox}

\title[Trigonal Morsifications on Hirzebruch surfaces Appendix by E. Shustin]{Trigonal Morsifications on Hirzebruch surfaces 
with an appendix by E. Shustin}
\author{Andrés Jaramillo Puentes}

\begin{document}
\input{morabstract}
\maketitle
\tableofcontents
\input{morintro}
\input{mor1}
\input{mor2}
\input{mor2c}
\input{appendix2}

\AtEndEnvironment{thebibliography}{

\bibitem{GLS} G.-M. Greuel, C. Lossen, and E. Shustin.
{\it Introduction to singularities and deformations}. Springer,
Berlin, 2007.

\bibitem{Hir} A. Hirschowitz. Le methode d'Horace pour l'interpolation \`a
plusieurs variables. {\it Manuscr. Math.} {\bf 50} (1985), 337--388.

\bibitem{KS} A. B. Korchagin and E. I. Shustin.
Affine curves of degree 6 and smoothings of ordinary 6th order singular point. {\it Math. USSR Izvestiya} {\bf 33}
(1989), no. 3, 501--520.

\bibitem{M} J. Milnor. \textit{Singular points of complex hypersurfaces}.
Princeton Univ. Press, Princeton, 1968.

\bibitem{OS} S. Yu. Orevkov and E. I. Shustin.
Real algebraic and pseudoholomorphic curves on the quadratic cone and smoothings of singularity $X_{21}$.
{\it St. Petersburg Math. J.} {\bf 28} (2017), 225--257.

\bibitem{Sh1} E. Shustin. Versal deformations in the space of plane curves of fixed degree.
{\it Function. Anal. Appl.} {\bf 21} (1987), 82--84.

\bibitem{ShT} E. Shustin. Gluing of singular and critical points.
\textit{Topology} {\bf 37} (1998), no. 1, 195--217.

\bibitem{Sh2} E. Shustin. A tropical approach to enumerative geometry.
\textit{St. Petersburg Math. J.} {\bf 17} (2006), 343--375.

\bibitem{ST} E. Shustin and I. Tyomkin. Patchworking singular algebraic curves, II.
\textit{Israel J. Math.} {\bf 151} (2006), 145--166.

\bibitem{AC1} N. A'Campo.  Le groupe de monodromie des
singularités isolées des courbes planes, I. Math. Ann. 213 (1975), 1--32.

\bibitem{AC2} N. A'Campo.  Real deformations and complex topology of plane 
curve singularities. Annales de la Faculté de Science de Toulouse, 6-e 
Ser., 8 (1999), no. 1, 5--23.

\bibitem{Ca} J. Callahan.  Singularities and Plane Maps. II. Sketching
Catastrofes. Amer. Math. Monthly 84 (1977), 765--803.

\bibitem{GZ1} S. M. Gusein-Zade.  Intersection matrices for certain 
singularities of functions
of two variables. Func. Anal. Appl. 8 (1974), no. 1, 10--13.

\bibitem{GZ2} S. M. Gusein-Zade.  Dynkin diagrams for singularities of 
functions of two variables.
Functional. Anal. Appl. 8 (1974), 295--300.

\bibitem{FPS} S. Fomin, P. Pylyavskyy, and E. Shustin.
Morsifications and mutations. Preprint at arXiv:1711.10598.

\bibitem{LS} P. Leviant and E. Shustin.  Morsifications of real
plane curve singularities. Preprint at arXiv:1703.05510.
}

\bibliographystyle{plain} 
\bibliography{bibliomor.bib} 
\end{document}

%% file: head.tex
\usepackage[T1]{fontenc}
\usepackage[utf8]{inputenc}
\usepackage[english]{babel}
\usepackage{amssymb,amsmath,amsthm,amscd}

\usepackage{caption}
\usepackage{subcaption}
\usepackage{url}
\usepackage[all]{xy}

\newcommand{\CP}{\mathbb{C}\mathbb{P}^{1}}
\newcommand{\RP}{\mathbb{R}\mathbb{P}^{1}}
\newcommand{\RPP}{\mathbb{R}\mathbb{P}^{2}}

\newcommand{\R}{\mathbb{R}}
\newcommand{\C}{\mathbb{C}}
\newcommand{\Z}{\mathbb{Z}}

\newcommand{\N}{\mathbb{N}}

\newcommand{\lra}{\longrightarrow}
\newcommand{\lmt}{\longmapsto}

\newcommand{\Si}{\Sigma}

\newtheorem{thm}{Theorem}[section]
\newtheorem{lm}[thm]{Lemma}
\newtheorem{coro}[thm]{Corollary}
\newtheorem{prop}[thm]{Proposition}
\theoremstyle{definition}
\newtheorem{df}[thm]{Definition}
\theoremstyle{plain}
\newtheorem{rmk}[thm]{Remark}

\newcommand{\tvs}{$\times$-vertices}
\newcommand{\gc}{generalized cut}

\DeclareMathOperator{\Cut}{Cut}
\DeclareMathOperator{\Ver}{Ver}
\DeclareMathOperator{\Ind}{Ind}
\DeclareMathOperator{\Dssn}{Dssn}

\DeclareMathOperator{\Sing}{Sing}

\usepackage{graphicx}
\usepackage{hyperref}\usepackage{tikz}
\usetikzlibrary{arrows,decorations.pathmorphing,decorations.pathreplacing,backgrounds,positioning,fit,petri,calc,cd}

%% file: morabstract.tex
\begin{abstract}
{
In this paper we obtain a classification of rigid isotopy classes of totally reducible trigonal curves lying on a Hirzebruch surface $\Sigma_n$, and having a maximal number of non-degenerated double points.
Such curves correspond to morsifications of a totally real semiquasihomogeneous singularity
of weight~$(3,3n)$ (the union of three smooth real branches intersecting each other with multiplicity~$n$).
We obtain this classification by studying combinatorial properties of dessins.

In the appendix, we prove that any morsification of a totally real semiquasihomogeneous singularity
of weight~$(3,3n)$ 
can be realized (up to isotopy) by the restriction of the equation to the Newton diagram and adding
monomials under the Newton diagram.}
\end{abstract}

%% file: morintro.tex
A \emph{morsification} of a real plane singularity is a real deformation with the maximal possible number of hyperbolic nodes.
Morsifications were introduced by N. A'Campo \cite{AC1,AC2} 
and S. Gusein-Zade \cite{GZ1, GZ2} as an important tool for the study of  
Dynkin diagrams, monodromy, topology of the singularity link, and other characteristics of singularities. Many interesting questions related  
to the geometry of morsifications still remain open, e.g., the problem of the existence of morsifications for arbitrary singularities (see \cite{LS}), or the relation  between different morsifications of singularities of the same topological type and mutational equivalence of quivers (see \cite{FPS}).
The present paper addresses the problem of isotopy classification of morfisications of singularities and completely solves the problem for singularities combined of three smooth real branches. 
Note that for simple ADE singularities the classification of morsifications is well known (see, for example, \cite{Ca}).

In Section~\ref{st:tri} we introduce dessins, a combinatorial tool encoding the geometry of trigonal curves lying on ruled surfaces.
We restrain ourselves to the case of real trigonal curves lying on the Hirzebruch surfaces~$\Sigma_n$ having~$3$ reducible real components and~$3n$ non-degenerated double points.
We refer to such curves as \emph{trigonal morsifications}. They represent morsifications of quasihomogeneous singularities of weight $(3,3n)$ (the union of three smooth real branches intersecting each other with multiplicity $n$) in the real plane.
We state combinatorial properties of the dessins associated to trigonal morsifications.

In Section~\ref{st:tm} we introduce \emph{wiring diagrams}, a sequential representation of the projective equivalence class of a trigonal morsification.
We use wiring diagrams to state Theorem~\ref{thm:wir}, which allow us to describe, in an inductive (on~$n$) manner, all projective equivalence classes of trigonal morsifications.
This description leads to the isotopy classification of trigonal morsifications.

Aditionally, we proved unicity for the rigid isotopy class of trigonal morsifications on $\Sigma_n$ whose wiring diagram has a maximal number of consecutive equal entries.
We introduce Reidemeister moves as a topological operation on the real part of a morsification and we prove that such operation can be realized algebraically, and that for every~$n$, there exists a unique isotopy class of trigonal morsification up to Reidemeister moves.

The Appendix~\ref{app}, by E. Shustin, is devoted to the semiquasihomogeneous case. 
We prove that any morsification of a totally real semiquasihomogeneous singularity
of weight $(3,3n)$ (the union of three smooth real branches intersecting each other with multiplicity $n$)
can be realized (up to isotopy) by a polynomial consisting of the restriction of the equation to the Newton diagram and adding
monomials under the Newton diagram.

{\small
{\bf Acknowledgements} The author would like to thank E. Shustin for the helpful discussions and for writing the Appendix.
The author's fellowship at the School of Mathematical Sciences of Tel Aviv University was supported by the grant no. 176/15 from the Israeli Science Foundation. 

The author would also like to thank the Institute Mittag-Leffler, Stockholm, for its support to visit during the research program Tropical Geometry, Amoebas and Polytopes.}

%% file: mor1.tex
\section{Trigonal curves and dessins}\label{st:tri}

In this section we introduce trigonal curves and dessins, which are the principal tool we use in order to study trigonal morsifications on the Hirzebruch surfaces $\Sigma_n$.
The content of this chapter is based on the book~\cite{deg} and the article~\cite{DIZ}.

\subsection{Ruled surfaces and trigonal curves}

\subsubsection{Basic definitions}
A compact complex surface $\Sigma$ is a \emph{(geometrically) ruled surface} over a curve $B$ 
if $\Si$ is endowed with a projection $\pi:\Si\lra B$ of fiber~$\CP$ as well as a special section $E$ 
of non-positive self-intersection.

\begin{df}
A \emph{reducible trigonal curve} is a curve $C$, lying in a ruled surface $\Sigma$ such that $C$ contains neither the exceptional section $E$ nor a fiber as component, and	the restriction $\pi|_{C}:C\lra B$ is a degree $3$ map.

A trigonal curve $C \subset \Sigma$ is {\it proper} if it does not intersect the exceptional section~$E$.
A {\it singular fiber} of a trigonal curve $C\subset\Si$ is a fiber~$F$ of $\Si$ intersecting $C\cup E$ geometrically in less than $4$ points.
\end{df} 


\subsubsection{Deformations} 
We are interested in the study of trigonal curves up to deformation. In the real case, we consider the curves up to equivariant deformation (with respect to the action of the complex conjugation, cf.~\ref{sssec:rtc}).

In the Kodaira-Spencer sense, a \emph{deformation} of the quintuple $(\pi\colon\Si\lra B,E,C)$ refers to an analytic space $X\lra S$ fibered over an marked open disk $S\ni o$ endowed with analytic subspaces $\mathcal{B}, \mathcal{E}, \mathcal{C}\subset X$ such that for every $s\in S$, the fiber~$X_{s}$ is diffeomorphic to $\Si$ and the intersections $\mathcal{B}_s:= X_s\cap \mathcal{B}$, $\mathcal{E}_s:= X_s\cap \mathcal{E}$ and $\mathcal{C}_s:= X_s\cap \mathcal{C}$ are diffeomorphic to $B$, $E$ and $C$, respectively, and there exists a map $\pi_s\colon X_s\lra B_s$ making $X_s$ a geometrically ruled surface over $B_s$ with exceptional section $E_s$, such that the diagram in Figure~\ref{fig:KodSpe} commutes
and $(\pi_o\colon X_o\lra B_o,E_o,C_o)=(\pi\colon\Si\lra B,E,C)$.

\begin{figure}[h]
\begin{center}
\begin{tikzcd}
	& E_s \arrow[r, "\text{diff.}"] \arrow[d, hook] & E \arrow[d, hook] & \\
C_s \arrow[rd, "\pi_s|_{C_s}"'] \arrow[r, hook]	& X_s \arrow[r, "\text{diff.}"] \arrow[d, "\pi_s"] & \Sigma	\arrow[d,"\pi"'] \arrow[r,hookleftarrow] & C \arrow[ld, "\pi|_C"]  \\
	& B_s \arrow[r, "\text{diff.}"] & B &
\end{tikzcd}
\end{center}
\caption{Commuting diagram of morphisms and diffeomorphisms of the fibers of a deformation.}
\label{fig:KodSpe}
\end{figure}

\begin{df}
 An \emph{elementary deformation} of a trigonal curve $C\subset\Si\lra B$ is a deformation of the quintuple $(\pi\colon\Si\lra B,E,C)$ in the Kodaira-Spencer sense.
\end{df}

An elementary deformation $X\lra S$ is \emph{equisingular} if for every $s\in S$ there exists a neighborhood $U_s\subset S$ of $s$ such that for every singular fiber~$F$ of $C$, there exists a neighborhood $V_{\pi(F)}\subset B$ of $\pi(F)$, where $\pi(F)$ is the only point with a singular fiber for every $t\in U_s$.
An elementary deformation over $D^2$ is a \emph{degeneration} or \emph{perturbation} if the restriction to $D^2\setminus\{0\}$ is equisingular and for a set of singular fibers~$F_i$ there exists a neighborhood $V_{\pi(F_i)}\subset B$ where there are no points with a singular fiber for every $t\neq0$. In this case we say that $C_t$ \emph{degenerates} to $C_0$ or $C_0$ is \emph{perturbed} to $C_t$, for $t\neq0$.

\subsubsection{Weierstraß equations}
 For a trigonal curve, the Weierstraß equations are an algebraic tool which allows us to study the behavior of the trigonal curve with respect to the zero section and the exceptional one. They give rise to an auxiliary morphism of $j$-invariant type, which plays an intermediary role between trigonal curves and dessins. 
Let $C\subset\Si\lra B$ be a proper trigonal curve. Mapping a point~$b\in B$ of the base to the barycenter of the points in $C\cap F_B^0$ (weighted according to their multiplicity) defines a section $B\lra Z\subset\Si$ called the {\it zero section}; it is disjoint from the exceptional section $E$.\\

The surface $\Si$ can be seen as the projectivization of a rank $2$ vector bundle, which splits as a direct sum of two line bundles such that the zero section $Z$ corresponds to the projectivization of $\mathcal{Y}$, one of the terms of this decomposition. In this context, the trigonal curve $C$ can be described by a Weierstraß equation, which in suitable affine charts has the form
\begin{equation} \label{eq:Weiertrass}
x^3+g_{2}x+g_{3}=0,
\end{equation}
where $g_{2}$, $g_{3}$ are sections of $\mathcal{Y}^2$, $\mathcal{Y}^3$ respectively, and $x$ is an affine coordinate such that $Z=\{x=0\}$ and $E=\{x=\infty\}$. For this construction, we can identify~$\Si\setminus B$ with the total space of $\mathcal{Y}$ and take $x$ as a local trivialization of this bundle. Nonetheless, the sections $g_{2}$, $g_{3}$ are globally defined. The line bundle $\mathcal{Y}$ is determined by $C$. The sections $g_{2}$, $g_{3}$ are determined up to change of variable defined by 
\begin{equation*}
(g_{2}, g_{3})\lra(s^2g_{2},s^3 g_{3} ),\;s\in H^{0}(B,\mathcal{O}^{*}_{B}).
\end{equation*}
Hence, the singular fibers of the trigonal curve $C$ correspond to the points where the equation~\eqref{eq:Weiertrass} has multiple roots, i.e., the zeros of the discriminant section
\begin{equation}
 \Delta:=-4g_2^3-27g_3^2\in H^0(B,\mathcal{O}_B(\mathcal{Y}^6)).
\end{equation}
A Nagata transformation over a point $b\in B$ changes the line bundle $\mathcal{Y}$ to~${\mathcal{Y}\otimes\mathcal{O}_B(b)}$ and the sections $g_2$ and $g_3$ to $s^2 g_2$ and $s^3 g_3$, where $s\in H^0(B,\mathcal{O}_B)$ is any holomorphic function having a zero at $b$.

\begin{df} \label{df:gen}
 Let $C$ be a non-singular trigonal curve with Weierstraß model determined by the sections $g_2$ and $g_3$ as in \eqref{eq:Weiertrass}. The trigonal curve $C$ is {\it almost generic} if every singular fiber corresponds to a simple root of the determinant section $\Delta=-4g_2^3-27g_3^2$ which is not a root of $g_2$ nor of $g_3$. The trigonal curve $C$ is {\it generic} if it is almost generic and the sections $g_2$ and $g_3$ have only simple roots. 
\end{df}

\subsubsection{The $j$-invariant}
The $j$-invariant describes the relative position of four points in the complex projective line $\CP$. We describe some properties of the $j$-invariant in order to use them in the description of the dessins.

\begin{df}
 Let $z_1$, $z_2$, $z_3$, $z_4\in\CP$. The $j$-{\it invariant} of a set $\{z_1, z_2, z_3, z_4\}$ is given by
\begin{equation} \label{eq:jinvariant}
\displaystyle j(z_1, z_2, z_3, z_4)=\frac{4(\lambda^2-\lambda+1)^3}{27\lambda^2(\lambda-1)^2},
\end{equation}
where $\lambda$ is the {\it cross-ratio} of the quadruple $(z_1, z_2, z_3, z_4)$ defined as

 \begin{equation*}
\displaystyle\lambda(z_1, z_2, z_3, z_4)=\frac{z_1-z_3}{z_2-z_3} : \frac{z_1-z_4}{z_2-z_4}.
\end{equation*}
 \end{df}
 
The cross-ratio depends on the order of the points while the $j$-invariant does not. Since the cross-ratio $\lambda$ is invariant under Möbius transformations, so is the $j$-invariant. When two points $z_i$, $z_j$ coincide, the cross-ratio $\lambda$ equals either $0$, $1$ or~$\infty$, and the $j$-invariant equals $\infty$. 
 
Let us consider a polynomial $z^3+g_2z+g_3$. We define the $j$-\emph{invariant} $j(z_1, z_2, z_3)$ of its roots $z_1$, $z_2$, $z_3$ as $j(z_1, z_2, z_3,\infty)$. If $\Delta=-4g_2^3-27g_3^2$ is the discriminant of the polynomial, then
\begin{equation*}
\displaystyle j(z_1, z_2, z_3,\infty)=\frac{-4g_2^3}{\Delta}.
\end{equation*}

A subset $A$ of $\CP$ is real if $A$ is invariant under the complex conjugation. We say that $A$ has a nontrivial symmetry if there is a nontrivial permutation of its elements which extends to a linear map $z\lmt az+b$, $a\in\C^*$, $b\in\C$.

\begin{lm}[\cite{deg}]
The set $ \{z_1, z_2, z_3\}$ of roots of the polynomial $z^3+g_2z+g_3$ has a nontrivial symmetry if and only if its $j$-invariant equals $0$ (for an order 3 symmetry) or 1 (for an order 2 symmetry).
\end{lm}

\begin{prop}[\cite{deg}]
Assume that $j(z_1, z_2, z_3)\in\R$. Then, the following holds
\begin{itemize}
\item
The $j$-invariant $j(z_1, z_2, z_3)<1$ if and only if the points $z_1, z_2, z_3$ form an isosceles triangle. The special angle seen as a function of the $j$-invariant is a increasing monotone function. This angle tends to $0$ when $j$ tends to $-\infty$, equals $\frac{\pi}{3}$ at $j=0$ and tends to $\frac{\pi}{2}$ when $j$ approaches $1$. 

\item
The $j$-invariant $j(z_1, z_2, z_3)\geq1$ if and only if the points $z_1, z_2, z_3$ are collinear.The ratio between the lengths of the smallest segment and the longest segment $\overline{z_lz_k}$ seen as a function of the $j$-invariant is a decreasing monotone function. This ratio equals $1$ when $j$ equals $1$, and $0$ when $j$ approaches $\infty$.
\end{itemize}
\end{prop}

%

\subsubsection{The $j$-invariant of a trigonal curve}
Let $C$ be a proper trigonal curve. We use the $j$-invariant defined for triples of complex numbers in order to define a meromorphic map $j_C$ on the base curve $B$. The map $j_{C}$ encodes the topology of the trigonal curve $C$. The map $j_C$ is called the $j$-\emph{invariant} of the curve $C$ and provides a correspondence between trigonal curves and dessins.

\begin{df}
 For a proper trigonal curve~$C$, we define its $j$-\emph{invariant} $j_C$ as the analytic continuation of the map
 \[
\begin{array}{ccc}
 B^{\#}&\lra&\C \\
 b &\lmt &j\mbox{-invariant of } C\cap F_b^0\subset F_b^0\cong\C.
\end{array}
\]
We call the trigonal curve $C$ {\it isotrivial} if its $j$-invariant is constant.
\end{df}

If a proper trigonal curve $C$ is given by a Weierstraß equation of the form \eqref{eq:Weiertrass}, then
\begin{equation}
j_C=-\frac{4g_2^3}{\Delta}\; , \text{ where }\Delta=-4g_2^3-27g_3^2.
\end{equation}

\begin{thm}[\cite{deg}] \label{thm:Etc}
 Let $B$ be a compact curve and $j\colon B\lra\CP$ a non-constant meromorphic map. Up to Nagata equivalence, there exists a unique trigonal curve $C\subset\Si\lra B$ such that $j_{C}=j$.
\end{thm}

Following the proof of the theorem, $j_B\lra\CP$ leads to a unique \emph{minimal} proper trigonal curve $C_j$, in the sense that any other trigonal curve with the same $j$-invariant can be obtained by positive Nagata transformations from $C_j$.

An equisingular deformation $C_s$, $s\in S$, of $C$ leads to an analytic deformation of the couple $(B_s,j_{C_s})$.

\begin{coro}[\cite{deg}] \label{coro:Etc}
 Let $(B,j)$ be a couple, where $B$ is a compact curve and $j\colon B\lra\CP$ is a non-constant meromorphic map. Then, any deformation of~$(B,j)$ results in a deformation of the minimal curve $C_j\subset\Si\lra B$ associated to $j$.
\end{coro}

The $j$-invariant of a generic trigonal curve $C\subset\Sigma\lra B$ has degree $\deg(j_C)=6d$, where $d=-E^2$. A positive Nagata transformation increases $d$ by one while leaving $j_C$ invariant.
The $j$-invariant of a generic trigonal curve $C$ has a ramification index equal to $3$, $2$ or $1$ at every point $b\in B$ such that $j_C(b)$ equals $0$, $1$ or $\infty$, respectively. We can assume, up to perturbation, that every critical value of $j_C$ is simple. In this case we say that $j_C$ has a {\it generic branching behavior}.

\subsubsection{Real trigonal curves}
\label{sssec:rtc}

We are mostly interested in real trigonal curves. A {\it real structure} on a complex variety $X$ is an anti-holomorphic involution $c\colon X\lra X$. We define a {\it real variety} as a couple $(X,c)$, where $c$ is a real structure on a complex variety $X$. We denote by $X_{\R}$ the fixed point set of the involution $c$ and we call $X_{\R}$ {\it the set of real points} of $c$.

A geometrically ruled surface $\pi\colon\Si\lra B$ is \emph{real} if there exist real structures $c_{\Si}\colon\Si\lra\Si$ and $c_{B}\colon B\lra B$ compatible with the projection $\pi$, i.e., such that $\pi\circ c_{\Si}=c_{B}\circ\pi$. We assume the exceptional section is {\it real} in the sense that it is invariant by conjugation, i.e., $c_{\Si}(E)=E$. 
Put $\pi_{\R}:=\pi|_{\Si_{\R}}\colon\Si_{\R}\lra B_{\R}$. Since the exceptional section is real, the fixed point set of every fiber is not empty, implying that the real structure on the fiber is isomorphic to the standard complex conjugation on $\CP$. Hence all the fibers of $\pi_{\R}$ are isomorphic to~$\RP$. 
Thus, the map $\pi_{\R}$ establishes a bijection between the connected components of the real part~$\Si_{\R}$ of the surface $\Si$ and the connected components of the real part~$B_{\R}$ of the curve~$B$. Every connected component of $\Si_{\R}$ is homeomorphic either to a torus or to a Klein bottle.

If $\Si=\mathbf{P}(1\oplus \mathcal{Y})$, with $\mathcal{Y}\in\operatorname{Pic}(B)$, we put $\mathcal{Y}_{i}:=\mathcal{Y}_{\R}|_{B_{i}}$ for every connected component $B_{i}$ of $B_{\R}$. Hence $\Si_{i}:=\Si_{\R}|_{B_{i}}$ is orientable if and only if $\mathcal{Y}_{i}$ is topologically trivial, i.e., its first Stiefel-Whitney class $w_{1}(\mathcal{Y}_{i})$ is zero.

\begin{df}
 A \emph{real trigonal curve} $C$ is a trigonal curve contained in a real ruled surface $(\Si,c_{\Si})\lra (B,c_{B})$ such that $C$ is $c_{\Si}$-invariant, i.e., $c_{\Si}(C)=C$.
\end{df}

If a real trigonal curve is proper, then $\mathcal{Y}$ is real as well as its $j$-invariant, seen as a morphism $j_{C}\colon(B,c_{B})\lra (\CP,z\lmt\bar{z})$, where $z\lmt\bar{z}$ denotes the standard complex conjugation on $\CP$. In addition, the sections $g_{2}$ and $g_{3}$ can be chosen real.

 Let us consider the restriction $\pi|_{C_{\R}}\colon C_{\R}\lra B_{\R}$. We put $C_{i}:=\pi|_{C_{\R}}^{-1}(B_{i})$ for every connected component $B_{i}$ of $B_{\R}$. We say that $B_{i}$ is \emph{hyperbolic} if $\pi|_{C_{i}}\colon C_{i}\lra~B_{i}$ has generically a fiber with three elements. The trigonal curve $C$ is \emph{hyperbolic} if its real part is non-empty and all the connected components of $B_{\R}$ are hyperbolic.
 
\subsection{Dessins}

The dessins d'enfants were introduced by A. Grothendieck (cf.~\cite{schneps})
in order to study the action of the absolute Galois group. We use a modified version of dessins d'enfants which was proposed by S. Orevkov~\cite{Orevkov}.

\subsubsection{Trichotomic graphs}

Let $S$ be a compact connected topological surface. A graph $D$ on the surface $S$ is a graph embedded into the surface and considered as a subset $D\subset S$. We denote by $\operatorname{Cut}(D)$ the \emph{cut} of $S$ along $D$, i.e., the disjoint union of the closure of connected components of $S\setminus D$.

\begin{df} \label{df:trigra}
 A \emph{trichotomic graph} on a compact surface $S$ is an embedded finite directed graph $D\subset S$ decorated with the following additional structures (referred to as \emph{colorings} of the edges and vertices of $D$, respectively):
 
\begin{itemize}
\item every edge of $D$ is color solid, bold or dotted,
\item every vertex of $D$ is black ($\bullet$), white ($\circ$), cross ($\times$) 
or monochrome (the vertices of the first three types are called \emph{essential}),
\end{itemize}
and satisfying the following conditions:
\begin{enumerate}
\item $\partial S\subset D$,
\item every essential vertex is incident to at least $2$ edges,
\item every monochrome vertex is incident to at least $3$ edges,
\item the orientations of the edges of $D$ form an orientation of the boundary $\partial\operatorname{Cut} (D)$ which is compatible with an orientation on $\operatorname{Cut} (D)$,
\item all edges incident to a monochrome vertex are of the same color,
\item $\times$-vertices are incident to incoming dotted edges and outgoing solid edges,
\item $\bullet$-vertices are incident to incoming solid edges and outgoing bold edges,
\item $\circ$-vertices are incident to incoming bold edges and outgoing dotted edges.
\end{enumerate}
\end{df}

Let $D\subset S$ be a trichotomic graph. A \emph{region} $R$ is an element of of $\operatorname{Cut(D)}$. The boundary $\partial R$ of $R$ contains $n=3k$ essential vertices. 
A region with $n$ essential vertices on its boundary is called an \emph{$n$-gonal region}.
We denote by $D_{solid}$, $D_{bold}$, $D_{dotted}$ the monochrome parts of $D$, i.e., the sets of vertices and edges of the specific color. On the set of vertices of a specific color, we define the relation $u\preceq v$ if there is a monochrome path from $u$ to~$v$, i.e., a path formed entirely of edge of the same color. We call the graph $D$ \emph{admissible} if the relation $\preceq$ is a partial order, equivalently, if there are no directed monochrome cycles.

\begin{df}
 A trichotomic graph $D$ is a \emph{dessin} if
\begin{enumerate}
 \item $D$ is admissible;
 \item every trigonal region of $D$ is homeomorphic to a disk.
\end{enumerate}
\end{df}

The orientation of the graph $D$ is determined by the pattern of colors of the vertices on the boundary of every region.

\subsubsection{Complex and real dessins}

Let $S$ be an orientable surface. Every orientation of $S$ induces a \emph{chessboard coloring} of $\Cut(D)$, i.e., a function on $\Cut(D)$ determining if a region $R$ endowed with the orientation set by $D$ coincides with the orientation of $S$.

\begin{df} A \emph{real trichotomic graph} on a real compact surface $(S,c)$ is a trichotomic graph $D$ on $S$ which is invariant under the action of $c$. Explicitly, every vertex $v$ of $D$ has as image $c(v)$ a vertex of the same color; every edge $e$ of $D$ has as image $c(e)$ an edge of the same color.
\end{df}

Let $D$ be a real trichotomic graph on $(S,c)$. Put $\overline{S}:=S/c$ as the quotient surface and put $\overline{D}\subset\overline{S}$ as the image of $D$ by the quotient map $S\lra S/c$. The graph $\overline{D}$ is a well defined graph on the surface $S/c$.

In the inverse sense, let $S$ be a compact surface, which can be non-orientable or can have non-empty boundary. Let $D\subset S$ be a trichotomic graph on $S$.
Consider its complex double covering $\widetilde{S}\lra S$ (cf. \cite{AG} for details), which has a real structure given by the deck transformation, and put $\widetilde{D}\subset\widetilde{S}$ the inverse image of $D$. The graph~$\widetilde{D}$ is a graph on $\widetilde{S}$ invariant by the deck transformation. We use these constructions in order to identify real trichotomic graphs on real surfaces with their images on the quotient surface.


\begin{prop}[\cite{deg}] \label{prop:realdd}
Let $S$ be a compact surface. Given a trichotomic graph ${D\subset S}$, then its oriented double covering $\widetilde{D}\subset\widetilde{S}$ is a real trichotomic graph. Moreover, $\widetilde{D}\subset\widetilde{S}$ is a dessin if and only if so is $D\subset S$.
Conversely, if $(S,c)$ is a real compact surface and $D\subset S$ is a real trichotomic graph, then its image $\overline{D}$ in the quotient $\overline{S}:=S/c$ is a trichotomic graph. Moreover, $\overline{D}\subset\overline{S}$ is a dessin if and only if so is $D\subset S$. 
\end{prop}

\begin{df}
Let $D$ be a dessin on a compact surface $S$. Let us denote by $\Ver(D)$ the set of vertices of $D$. For a vertex $v\in\Ver(D)$, we define the \emph{index} $\Ind(v)$ of $v$ as half of the number of incident edges of $\widetilde{v}$, where $\widetilde{v}$ is a preimage of $v$ by the double complex cover of $S$ as in Proposition~\ref{prop:realdd}.

A vertex $v\in\Ver(D)$ is \emph{singular} if
\begin{itemize}
\item $v$ is black and $\Ind(v)\not\equiv0\mod3$,
\item or $v$ is white and $\Ind(v)\not\equiv0\mod2$,
\item or $v$ has color $\times$ and $\Ind(v)\geq2$.
\end{itemize}
We denote by $\Sing(D)$ the set of singular vertices of $D$. A dessin is \emph{non-singular} if none of its vertices is singular.
\end{df}

\begin{df}
 Let $B$ be a complex curve and let $j\colon B\lra\CP$ a non-constant meromorphic function, in other words, a ramified covering of the complex projective line. The dessin $D:=\Dssn(j)$ associated to $j$ is the graph given by the following construction:
\begin{itemize}
\item as a set, the dessin $D$ coincides with $j^{-1}(\RP)$, where $\RP$ is the fixed point set of the standard complex conjugation in $\CP$;
\item black vertices $(\bullet)$ are the inverse images of $0$;
\item white vertices $(\circ)$ are the inverse images of $1$;
\item vertices of color $\times$ are the inverse images of $\infty$;
\item monochrome vertices are the critical points of $j$ with critical value in $\CP\setminus\{0,1,\infty\}$;
\item solid edges are the inverse images of the interval $[\infty,0]$;
\item bold edges are the inverse images of the interval $[0,1]$;
\item dotted edges are the inverse images of the interval $[1,\infty]$;
\item orientation on edges is induced from an orientation of $\RP$.
\end{itemize}
\end{df}

\begin{lm}[\cite{deg}]
Let $S$ be an oriented connected closed surface. Let ${j\colon S\lra\CP}$ a ramified covering map. The trichotomic graph $D=\Dssn(j)\subset S$ is a dessin. Moreover, if $j$ is real with respect to an orientation-reversing involution $c\colon S\lra S$, then $D$ is $c$-invariant.
\end{lm}

Let $(S,c)$ be a compact real surface. If $j\colon (S,c)\lra (\CP,z\lra\bar{z})$ is a real map, we define $\Dssn_{c}(j):=\Dssn(j)/c\subset S/c$.

\begin{thm}[\cite{deg}]
Let $S$ be an oriented connected closed surface (and let \linebreak[4] ${c\colon S\lra S}$ an orientation-reversing involution). A (real) trichotomic graph ${D\subset S}$ is a (real) dessin if and only if $D=\Dssn(j)$ 
for a (real) ramified covering \linebreak[4]  ${j\colon S\lra\CP}$. 

Moreover, $j$ is unique up to homotopy in the class of (real) ramified coverings with dessin $D$.
\end{thm}

The last theorem together with the Riemann existence theorem provides the next corollaries, for the complex and real settings.

\begin{coro}[\cite{deg}] \label{coro:Ej}
Let $D\subset S$ be a dessin on a compact closed orientable surface~$S$. Then there exists a complex structure on $S$ and a holomorphic map ${j\colon S\lra\CP}$ such that $\Dssn(j)=D$. Moreover, this structure is unique up to deformation of the complex structure on $S$ and the map $j$ in the Kodaira-Spencer sense.
\end{coro}

\begin{coro}[\cite{deg}] \label{coro:EjR}
Let $D\subset S$ be a dessin on a compact surface $S$. Then there exists a complex structure on its double cover $\widetilde{S}$ and a holomorphic map $j\colon\widetilde{S}\lra\CP$ such that $j$ is real with respect to the real structure $c$ of $\widetilde{S}$ and $\Dssn_{c}(j)=D$. Moreover, this structure is unique up to equivariant deformation of the complex structure on $S$ and the map $j$ in the Kodaira-Spencer sense.
\end{coro}

\subsubsection{Deformations of dessins}

In this section we describe the notions of deformations which allow us to associate classes of non-isotrivial trigonal curves and classes of dessins, up to deformations and equivalences that we explicit.

\begin{df}
A {\it deformation of coverings} is a homotopy $S\times [0,1] \lra \CP$ within the class of (equivariant) ramified coverings. The deformation is {\it simple} if it preserves the multiplicity of the inverse images of $0$, $1$, $\infty$ and of the other real critical values.
\end{df}

Any deformation is locally simple except for a finite number of values~$t\in [0,1]$.

\begin{prop}[\cite{deg}]
 Let $j_{0}, j_{1}\colon S\lra\CP$ be ($c$-equivariant) ramified coverings. They can be connected by a simple (equivariant) deformation if and only their dessins $D(j_0)$ and $D(j_{1})$ are isotopic (respectively, $D_{c}(j_{0})$ and $D_{c}(j_{1})$).
\end{prop}

\begin{df} \label{df:equidef}
 A deformation $j_{t}\colon S\lra\CP$ of ramified coverings is \emph{equisingular} if the union of the supports 
\[\bigcup_{t\in [0,1]} \operatorname{supp}\left\{ (j^{*}_{t}(0)\mod 3)+(j^{*}_{t}(1)\mod 2)+j^{*}_{t}(\infty)\right\}\] 
considered as a subset of $S\times [0,1]$ is an isotopy. Here $^{*}$ denotes the divisorial pullback of a map $\varphi:S\lra S'$ at a point $s'\in S'$:
 \[ \varphi^{*}(s')=\displaystyle\sum_{s\in\varphi^{-1}(s')}r_{s}s, \] where $r_{s}$ if the ramification index of $\varphi$ at $s\in S$.
\end{df}

A dessin $D_1\subset S$ is called a \emph{perturbation} of a dessin $D_0\subset S$, and $D_0$ is called a \emph{degeneration} of $D_1$, if for every vertex $v\in \Ver(D_0)$ there exists a small neighboring disk $U_v\subset S$ such that $D_0\cap U_v$ only has edges incident to $v$ and $D_1\cap U_v$ contains essential vertices of at most one color. 

\begin{thm}[\cite{deg}] \label{thm:support}
 Let $D_{0}\subset S$ be a dessin, and let $D_{1}$ be an admissible perturbation. Then there exists a map $j_{t}\colon S\lra\CP$ such that
\begin{enumerate}
 \item $D_{0}=\Dssn(j_0)$ and $D_{1}=\Dssn(j_{1})$;
 \item $j_{t}|_{S\setminus\bigcup_{v}U_{v}}=j_{t'}|_{S\setminus\bigcup_{v}U_{v}}$ for every $t$, $t'$;
 \item the deformation restricted to $S\times (0,1]$ is simple.
\end{enumerate}
\end{thm}

\begin{coro}[\cite{deg}] \label{coro:def}
 Let $S$ be a complex compact curve, $j\colon S\lra\CP$ a non-constant holomorphic map, and let $\Dssn_{c}(j)=D_{0}, D_{1}, \dots, D_{n}$ be a chain of dessins in $S$ such that for $i=1, \dots, n$ either $D_{i}$ is a perturbation of $D_{i-1}$, or $D_{i}$ is a degeneration of $D_{i-1}$, or $D_{i}$ is isotopic to $D_{i-1}$. Then there exists a piecewise-analytic deformation $j_{t}\colon S_{t}\lra\CP$, $t\in[0,1]$, of $j_0=j$ such that $\Dssn(j_{1})=D_{n}$.
\end{coro}

\begin{coro}[\cite{deg}] \label{coro:defR}
 Let $(S,c)$ be a real compact curve, $j\colon (S,c)\lra(\CP,{z\lmt\bar{z}})$ be a real non-constant holomorphic map, and let $\Dssn_{c}(j)= D_{0}, D_{1}, \dots, D_{n}$ be a chain of real dessins in $(S,c)$ such that for $i=1, \dots, n$ either $D_{i}$ is a equivariant perturbation of $D_{i-1}$, or $D_{i}$ is a equivariant degeneration of $D_{i-1}$, or $D_{i}$ is equivariantly isotopic to $D_{i-1}$. Then there is a piecewise-analytic real deformation $j_{t}\colon (S_{t},c_{t})\lra(\CP,\bar{\cdot})$, $t\in[0,1]$, of $j_0=j$ such that $\Dssn_{c}(j_{1})=D_{n}$.
\end{coro}

Due to Theorem \ref{thm:support}, the deformation $j_t$ given by Corollaries \ref{coro:def} and \ref{coro:defR} is equisingular in the sense of Definition \ref{df:equidef} if and only if all perturbations and degenerations of the dessins on the chain $D_{0}, D_{1}, \dots, D_{n}$ are equisingular.

\subsection{Trigonal curves and their dessins}

In this section we describe an equivalence between dessins.

\subsubsection{Correspondence theorems}

Let $C\subset\Si\lra B$ be a non-isotrivial proper trigonal curve.
We associate to $C$ the dessin corresponding to its $j$-invariant $\Dssn(C):=\Dssn(j_{C})\subset B$. In the case when $C$ is a real trigonal curve we associate to $C$ the real dessin corresponding to its $j$-invariant, $\Dssn_{c}(C):=\Dssn(j_{C})\subset B/{c_{B}}$, where~$c_{B}$ is the real structure of the base curve $B$.\\

So far, we have focused on one direction of the correspondences: we start with a trigonal curve $C$, 
consider its $j$-invariant and 
construct the dessin associated to it. 
Now, we study the opposite direction. 
Let us consider a dessin $D$ on a topological orientable closed surface $S$. By Corollary \ref{coro:Ej},
there exist a complex structure~$B$ on~$S$ and a holomorphic map $j_{D}\colon B\lra\CP$ such that $\Dssn(j_{D})=D$. By Theorem~\ref{thm:Etc} and Corollary~\ref{coro:Etc} there exists a trigonal curve~$C$ having~$j_{D}$ as $j$-invariant; such a curve is unique up to deformation in the class of trigonal curves with fixed dessin. Moreover, due to Corollary~\ref{coro:def}, any sequence of isotopies, perturbations and degenerations of dessins gives rise to a piecewise-analytic deformation of tri\-gonal curves, which is singular if and only if all perturbations and degenerations are.

In the real framework, let $(S,c)$ a compact close oriented topological surface endowed with a orientation-reversing involution. Let $D$ be a real dessin on $(S,c)$. By Corollary~\ref{coro:EjR}, there exists a real structure $(B,c_{B})$ on $(S,c)$ and a real map ${j_{D}\colon(B,c_{B})\lra(\CP,z\longmapsto\bar{z})}$ such that $\Dssn_c(j_{D})=D$. By Theorem \ref{thm:Etc}, Corollary \ref{coro:Etc} and the remarks made in Section \ref{sssec:rtc}, there exists a real trigonal curve~$C$ having~$j_{D}$ as $j$-invariant; such a curve is unique up to equivariant deformation in the class of real trigonal curves with fixed dessin. Furthermore, due to Corollary~\ref{coro:defR}, any sequence of isotopies, perturbations and degenerations of dessins gives rise to a piecewise-analytic equivariant deformation of real trigonal curves, which is equisingular if and only if all perturbations and degenerations are.

\begin{df}
A dessin is \emph{reduced} if
\begin{itemize}
\item 
for every $v$ $\bullet$-vertex one has $\Ind{v}\leq3$,
\item 
for every $v$ $\circ$-vertex one has $\Ind{v}\leq2$,
\item every monochrome vertex is real and has index $2$.
\end{itemize}
A reduced dessin is {\it generic} if all its $\bullet$-vertices and $\circ$-vertices are non-singular and all its $\times$-vertices have index~$1$.
\end{df}

Any dessin admits an equisingular perturbation to a reduced dessin. The vertices with excessive index (i.e., index greater than 3 for $\bullet$-vertices or than 2 for $\circ$-vertices) can be reduced by introducing new vertices of the same color.

\begin{figure} 
\centering
\includegraphics[width=5in]{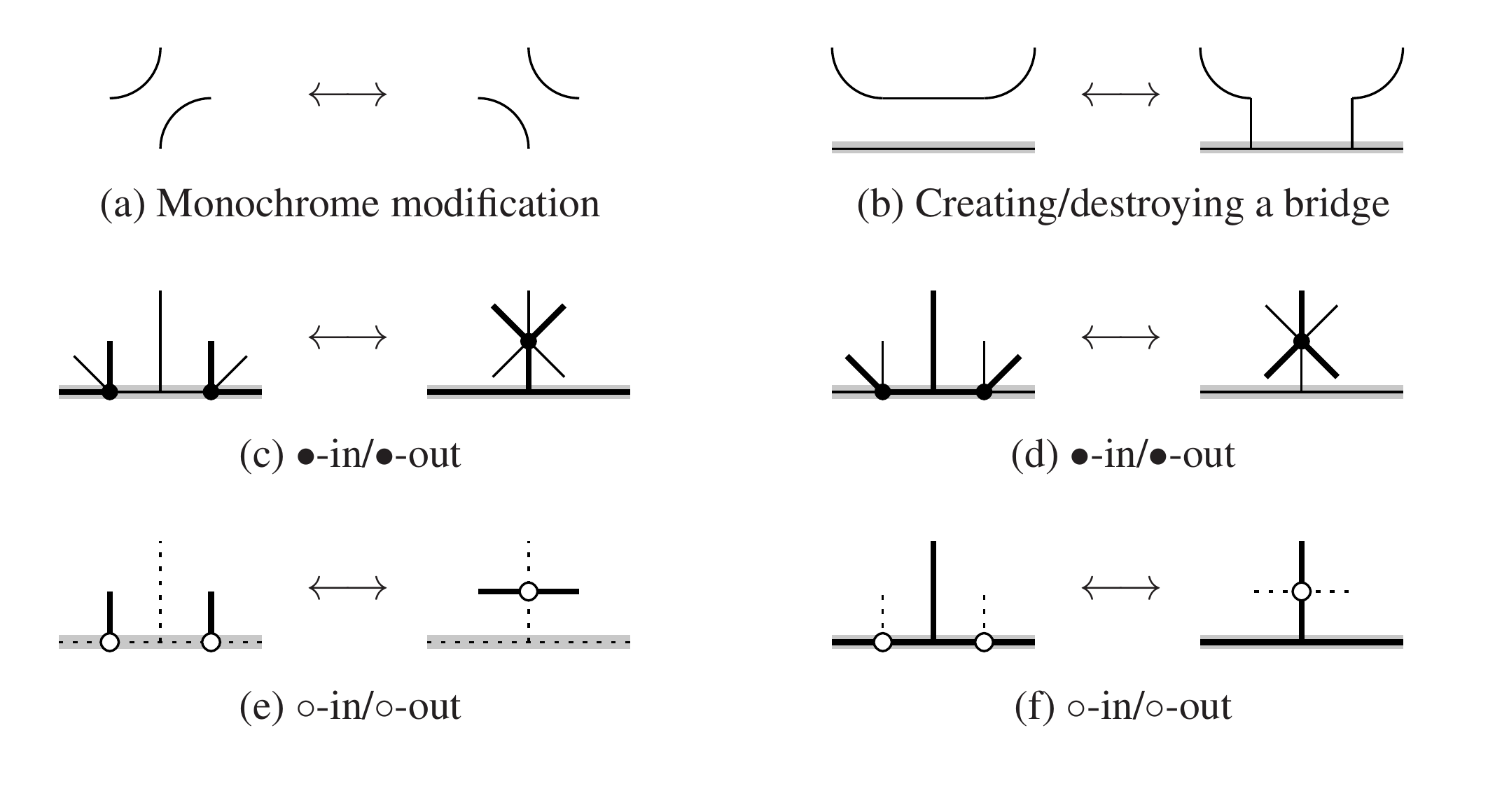}
\caption{Elementary moves.\label{fig:elem}}
\end{figure}

In order to define an equivalence relation of dessins, we introduce 
{\it elementary moves}. Consider two reduced dessins $D$, $D'\subset S$ such that they coincide outside a closed disk $V\subset S$.
If $V$ does not intersect $\partial S$ and the graphs  $D\cap V$ and $D'\cap V$ are as shown in Figure~\ref{fig:elem}(a), then we say that performing a {\it monochrome modification} on the edges intersecting $V$ produces $D'$ from $D$, or {\it vice versa}.
This is the first type of {\it elementary moves}. 
Otherwise, the boundary component inside $V$ is shown in light gray. In this setting, if the graphs $D\cap V$ and $D'\cap V$ are as shown in one of the subfigures in Figure~\ref{fig:elem}, we say that performing an \emph{elementary move} of the corresponding type on $D\cap V$ produces $D'$ from $D$, or {\it vice versa}. 

\begin{df} 
Two reduced dessins $D$, $D'\subset S$ are {\it elementary equivalent} if, after a (preserving orientation, in the complex case) homeomorphism of the underlying surface $S$ they can be connected by a sequence of isotopies and elementary moves between dessins, as described in Figure~\ref{fig:elem}.
\end{df}

This definition is meant so that two reduced dessins are elementary equivalent if and only if they can be 
connected up to homeomorphism by a sequence of isotopies, equisingular perturbations and degenerations.

The following theorems establish the equivalences between the deformation classes of trigonal curves we are interested in and elementary equivalence classes of certain dessins. We use these links to obtain different classifications of curves {\it via} the combinatorial study of dessins.

 \begin{thm}[\cite{deg}] \label{th:correpondance2}
 There is a one-to-one correspondence between the set of equivariant equisingular deformation classes of non-isotrivial proper real trigonal curves $C\subset\Si\lra(B,c)$ with $\widetilde{A}$ type singular fibers only and the set of elementary equivalence classes of reduced real dessins $D \subset B/c$.
\end{thm}

This correspondences can be extended to trigonal curves with more general singular fibers (see \cite{deg}).

\begin{df} 
Let $C\subset\Si\lra B$ be a proper trigonal curve. 
We define the {\it degree} of the curve $C$ as $\deg(C):=-3E^2$ where $E$ is the exceptional section of~$\Si$. For a dessin $D$, we define its {\it degree} as $\deg(D)=\deg(C)$ where $C$ is a minimal proper trigonal curve such that $\Dssn(C)=D$.
\end{df}

\subsubsection{Real generic curves}

Let $C$ be a generic real trigonal curve and let $D:=\Dssn_{c}(C)$ be a generic dessin. 
The {\it real part} of $D\subset S$ is the intersection $D\cap\partial S$. For a specific color $\ast\in\{\mathrm{ solid},\mathrm{ bold},\mathrm{ dotted}\}$, $D_{\ast}$ is the subgraph of the corresponding color and its adjacent vertices. The components of $D_{\ast}\cap\partial S$ are either components of~$\partial S$, called \emph{monochrome components} of $D$, or 
segments, called {\it maximal monochrome segments} of $D$. 
We call these monochrome components or segments \emph{even} or \emph{odd} according to 
the parity of the number of $\circ$-vertices they contain.

Moreover, we refer to the dotted monochrome components as \emph{hyperbolic components}.

\subsection{Toiles}\label{ch:toi}

Within the moduli space of trigonal curves of fixed degree, generic trigonal curves are smooth and the discriminant (i.e., the set of singular trigonal curves, cf. \cite{DK}) has a stratification which we can describe by means of dessins. 
Singular proper trigonal curves have singular dessins and the singular points are represented by singular vertices. 
A generic singular trigonal curve $C$ has exactly one singular point, which is a non-degenerate double point ({\it node}). 
Moreover, if $C$ is a proper trigonal curve, then the double point on it is represented by a $\times$-vertex of index $2$ on its dessin.
In addition, if $C$ has a real structure,
the double point is real and so is its corresponding vertex, leading to the cases where the $\times$-vertex of index $2$ has dotted real edges (representing the intersection of two real branches) or has solid real edges (representing one isolated real point, which is the intersection of two complex conjugated branches).

\begin{df} 
Let $D\subset S$ be a dessin on a compact surface $S$. A {\it nodal vertex} 
({\it node})  
of $D$ is a $\times$-vertex of index $2$.
The dessin $D$ is called \emph{nodal} if all its singular vertices are nodal vertices.
We call a \emph{toile} a non-hyperbolic real nodal dessin on $(\CP,z\lmt\bar{z})$. 
\end{df}

Since a real dessin on $(\CP,z\lmt\bar{z})$ descends to the quotient, we represent toiles on the disk.

In a real dessin, there are two types of real nodal vertices, namely, vertices having either real solid edges and interior dotted edges, or dotted real edges and interior solid edges. We call \emph{isolated nodes} of a dessin $D$ thoses $\times$-vertices of index~$2$ corresponding to the former case and \emph{non-isolated nodes} those corresponding to the latter.

\begin{df} 
Given a dessin $D$, a subgraph $\Gamma\subset D$ is a \emph{cut} if it consists of a single interior edge connecting two real monochrome vertices. An {\it axe} is an interior edge of a dessin connecting a $\times$-vertex of index $2$ and a real monochrome vertex.
\end{df}
 
Let us consider a dessin $D$ lying on a surface $S$ having a cut or an axe~$T$. 
Assume that $T$ divides $S$ and consider the connected components $S_{1}$ and $S_2$ of~$S\setminus T$. 
Then, we can define two dessins $D_1$, $D_2$, each lying on the compact surface $\overline{S_i}\subset S$, respectively for $i=1, 2$, and determined by $D_i:=(D\cap S_i)\cup\{T\}$. 
If $S\setminus T$ is connected, we define the surface $S'=(S\setminus T)\sqcup T_1 \sqcup T_2/\varphi_1, \varphi_2$, where $\varphi_i:T_i\lra S$ is the inclusion of one copy $T_i$ of $T$ into $S$, and the dessin $D':=(D\setminus T)\sqcup T_1 \sqcup T_2/\varphi_1, \varphi_2$.
 
By these means, a dessin having a cut or an axe determines either two other dessins of smaller degree or a dessin lying on a surface with a smaller fundamental group. 
Moreover, in the case of an axe, the resulting dessins have one singular vertex less. 
Considering the inverse process, we call $D$ the \emph{gluing} of $D_1$ and $D_2$ along $T$ or the \emph{gluing} of $D'$ with itself along $T_1$ and $T_2$. 

\begin{df} \label{df:gennodal}
Let $C\subset\Sigma$ 
be a nodal proper trigonal curve with nodes $n_1, n_2, \dots, n_l$. 
Consider a Weierstraß model of $C$ determined by 
sections $g_2$ and $g_3$.
We say that the trigonal curve $C$ is \emph{almost generic} if every singular fiber different from the fiber at $\pi(n_i)$ corresponds to a simple root of the determinant section $\Delta=-4g_2^3-27g_3^2$ which is not a root of $g_2$ nor of $g_3$. The nodal trigonal curve $C$ is \emph{generic} if it is almost generic and the sections $g_2$ and $g_3$ only have simple roots.
\end{df}

\label{sec:toiles}

\begin{df} 
Given a real dessin $D$ and a vertex $v\in\Ver(D)$, we call the \emph{depth} of~$v$ the minimal number $n$ such that there exists an undirected inner chain $v_0,\dots, v_n$ in $D$ from $v_0=v$ to a real vertex $v_n$ and we denote the depth of $v$ by $\operatorname{dp}(v)$.
The depth of a dessin $D$ is defined as the maximum of the depth of the black and white vertices of $D$ and it is denoted by $\operatorname{dp}(D)$.
\end{df}

\begin{df} 
A {\it {\gc}} of a dessin $D$ is an inner undirected chain formed entirely of inner edges of the same color, either dotted or solid, connecting two distinct real nodal or monochrome vertices.
\end{df}

Analogously to a cut, cutting a dessin $D\subset S$ by a {\gc} produces two dessins of lower degree or a dessin of the same degree in a surface with a simpler topology, depending on whether the inner chain divides or not the surface $S$.

\begin{prop}[\cite{deg5}] \label{prop:deep}
Let $D\subset \mathbb{D}^2$ be a toile of degree greater than $3$. Then, there exists a toile $D'$ weakly equivalent to $D$ such that either $D'$ has depth $1$ or $D'$ has a {\gc}.
\end{prop}

\begin{coro} \label{cr:deepbf}
Let $D$ be as in Proposition \ref{prop:deep}. If there exists a toile $D'$ weakly equivalent to $D$ with depth $1$, then $D'$ can be chosen bridge-free.
\end{coro}

\begin{prop}[\cite{deg5}] \label{prop:decomptoi}
Let $D$ be a toile of degree at least $6$ and depth at most $1$.
Then, there exists a toile $D'$ weakly equivalent to $D$ such that $D'$ has a {\gc}.
Moreover, if $D$ has isolated real nodal {\tvs}, the {\gc} is dotted or a solid axe.
\end{prop}

%% file: mor2.tex
\section{Trigonal morsifications on Hirzebruch surfaces}
\label{st:tm}

\begin{df}
A \emph{trigonal morsification} on $\Sigma_n$ is a hyperbolic trigonal curve on~$\Sigma_n$ having as singular points exactly $3n$ real nodal points.
\end{df}

Henceforth, we only consider trigonal morsifications on the Hirzebruch surfaces~$\Sigma_n$.

\subsection{Wiring diagrams}
Given a morsification $C$ and a real white vertex \linebreak[4] ${v\in\Dssn(C)}$, we associate to $v$ a sequence $(a_i)_{i=1}^{3n}$ in the following way.
First, let us fix an orientation of the real part of the dessin.
Then, enumerate accordingly the set of nodal {\tvs} as $\{u_1, u_2,\dots, u_{3n}\}$.
Lastly, put the number $a_i=1$ if the arc~$(v,u_i)$ has an even number of white vertices, and $a_i=2$ otherwise.

Since the elementary moves of dessins do not modify the parity of segment between nodal {\tvs}, this sequence is well defined in the equivalence class of dessins.

If $n$ is odd, the sequence $(a_{3n-i})_{i=1}^{3n}$ corresponds to the opposite orientation.
On the other hand, if $n$ is even, the sequence $(a_i')_{i=1}^{3n}$, with $a_i'=3-a_{3n-i}$, corresponds to the opposite orientation.
We say that two sequences $(a_i)_{i=1}^{3n}, (b_i)_{i=1}^{3n}$ are \emph{equivalent} if they are equal, if $b_i=3-a_i, i=1,2,\dots,3n$, if $b_i=a_{3n-i}, i=1,2,\dots,3n$ or if~$b_i=3-a_{3n-i}, i=1,2,\dots,3n$.

We call the \emph{wiring diagram} of $D$ with respect to $v$ the sequence $(a_i)_{i=1}^{3n}$ associated to $v$ up to equivalence.

We put $\sigma_m\in\mathcal{S}_{m}$ as the permutation
\[\sigma_m=\left(\begin{array}{ccccc}
1 & 2 & \dots & m-1 & m\\
2 & 3 & \dots & m   & 1
\end{array} \right),\]

and we define $\sigma_{3n}\cdot(a_i)_{i=1}^{3n}$ as the sequence given by $(a_{\sigma_{3n}(i)})_{i=1}^{3n}$ for $n$  even, or by the sequence
\[(\sigma_{3n}\cdot(a_i)_{i=1}^{3n})_j=
\left\{\begin{array}{ll}
a_{\sigma_{3n}(j)} 	& \text{ if $j\neq 3n$,}\\
3-a_{1}	& \text{ if $j=3n$,}
\end{array} \right.
\]
for $n$ odd.

We say that two sequences $(a_i)_{i=1}^{3n}, (b_i)_{i=1}^{3n}$ are \emph{projective equivalent} if there exists~$k\in\N$ such that $(b_i)_{i=1}^{3n}$ is equivalent to $\sigma_{3n}^k\cdot (a_i)_{i=1}^{3n}$.

\begin{lm}
Let $C$ be a trigonal morsification and let $D=\Dssn(C)$ be its dessin.
If~$v$, $v'$ are real white vertices of $D$, then, the wiring diagrams of $D$ with respect to~$v$ and~$v'$ are projective equivalent.
\end{lm}

\begin{proof}
Fix an orientation of $\RPP$, the boundary of $D$.
Let $k$ and $l$ be the number of nodal {\tvs} and white vertices, respectively, in the arc $(v,v']$.\linebreak
Let~$(a_i)_{i=1}^{3n}, (a_i ')_{i=1}^{3n}$ the wiring diagrams of $D$ with respect to $v$ and $v'$, respectively.
If~$n$ is even, then $a_j '=a_j$ if $l$ is even or $a_j '=3-a_j$ if $l$ is odd, for $1\leq j\leq 3n$.
If~$n$ is odd, then $a_j '=a_j$, $1\leq j\leq 3n-k$ and $a_j '=3-a_j$, $3n-k<j$ if $l$ is even or~$a_j '=3-a_j$, $1\leq j\leq 3n-k$ and $a_j '=a_j$, $3n-k<j$ if $l$ is odd.
Therefore, we have that $(a_i ')_{i=1}^{3n}=\sigma_{3n}^k (a_i)_{i=1}^{3n}$ or $(a_i ')_{i=1}^{3n}=\sigma_{3n}^k (3-a_i)_{i=1}^{3n}$, depending whether $l$ is even or odd, respectively.
\end{proof}

\begin{df}
Given a morsification $C$ on $\Sigma_n$, we define as the \emph{wiring diagram} of $C$ the projective equivalence class of the wiring diagram of $\Dssn(D)$ with respect to any real white vertex.
\end{df}

Given a sequence $(a_i)_{i=1}^{3n}$, we define its \emph{length vector} as the vector $o$ whose $i$-th entry $o_i$ correspond to the length of the $i$-th maximal set of indexes \linebreak[4]
${I_i=\{i_0,{i_0+1},\dots,i_0+o_i\}}$ such that $a_j=a_l, \forall j,l\in I_i$.
The vector $o$ is an ordered partition of $3n$.
We denote by $o(a_i)$ the entry in $o$ correponding to the maximal set of indexes containing $i$.

Given a projective equivalence class of a sequence $(a_i)_{i=1}^{3n}$, we define its \emph{cyclic length vector} as the length vector $o$ of a projective equivalent sequence $(a_i')_{i=1}^{3n}$ such that $a_1\neq a_{3n}$ if $n$ is even, or such that $a_1=a_{3n}$ if $n$ is odd.
(If such a sequence does not exist, we put its cyclic length vector as $(3n)$).
We consider this vector $o$ as an element of $\Z^{m}/<\sigma_m,\tau_m>$, where $\tau_m\in\mathcal{S}_m$ is the permutation
\[\tau_m=\left(\begin{array}{ccccc}
1 & 2 & \dots & m-1 & m\\
m & m-1 & \dots & 2   & 1
\end{array} \right),\]
so in that manner, it does not depend on the representative of the projective equivalent class.

\begin{df}
Given a morsification $C$ on $\Sigma_n$, we define as the \emph{cyclic length vector} of $C$ the cyclic length vector of its wiring diagram.
\end{df}
 
\begin{lm}
Every projective equivalence class of sequences is determined by its cyclic length vector.
\end{lm}

\begin{proof}
Let $(a_i)_{i=1}^{3n}$, $(b_i)_{i=1}^{3n}$ be two sequences such that they have the same cyclic lenght vector $o$.
After a shifting $\sigma_{3n}^k, k\in\{0,1,...,3n-1\}$ and a permutation $\tau_{3n}^l, l\in\{0,1\}$, the sequences $(a_i)_{i=1}^{3n}$ and $\sigma_{3n}^k\circ\tau_{3n}^l\circ(b_i)_{i=1}^{3n}$ have the same length vector.
Then, either $a_i=\sigma_{3n}^k\circ\tau_{3n}^l(b_i)$,${i=1,2,\dots,3n}$ or $a_i=3-\sigma_{3n}^k\circ\tau_{3n}^l(b_i), {i=1,2,\dots,3n}$.
Thus, the sequences $(a_i)_{i=1}^{3n}$ and $(b_i)_{i=1}^{3n}$ are projective equivalent.
\end{proof}

\subsection{Decomposition of dessins and wiring diagrams}

\begin{figure}
\begin{subfigure}[b]{0.4\textwidth}
\includegraphics[width=\textwidth]{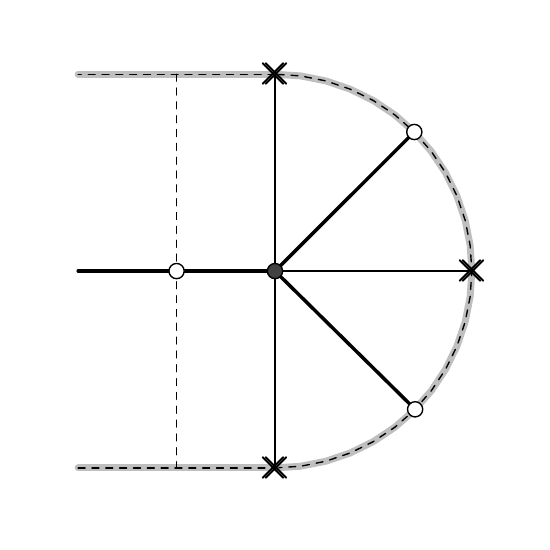}
\caption{Dotted cut resulting in a cubic $H^{***}$}
\label{fg:gludot}
\end{subfigure}
\begin{subfigure}[b]{0.4\textwidth}
\includegraphics[width=\textwidth]{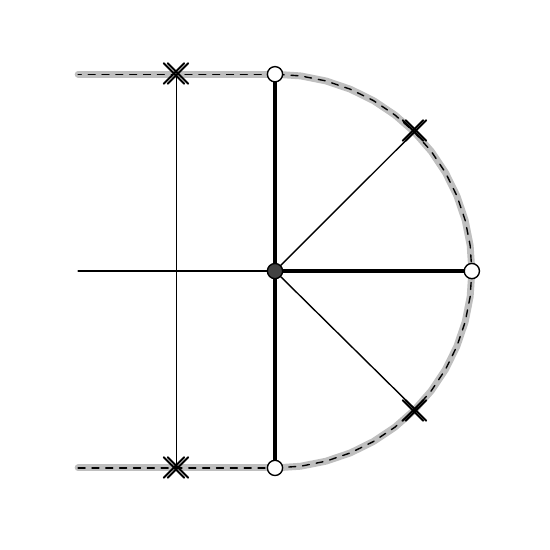}
\caption{Solid generalized cut resulting in a cubic $I_{2}^{**}$}
\label{fg:glusol}
\end{subfigure}
\caption{Block decomposition of dessins}
\end{figure}

Due to Proposition \ref{prop:decomptoi}, the dessin associated to a morfisication $C\subset\Sigma_n$ is elementary equivalent to the gluing of $n$ cubic dessins. 
There are two kinds of gluings: through a \emph{dotted cut} (see Figure \ref{fg:gludot}) or through a \emph{solid generalized cut} (see Figure \ref{fg:glusol}).

In the case the gluing is though a dotted cut, the dessin $D$ of a morsification is elementary equivalent to the gluing $D'-H^{***}$, where $H^{***}$ is the hyperbolic cubic dessin corresponding to three real lines in general position in $\RPP$, after a strict transform with respect to the blow up of a point in $\RPP$ which does not lie on any line.
Then the wiring diagram of $D$ would have either $1121$ or $2212$ corresponding to one of the nodal {\tvs} of $D'$ and the nodal {\tvs} of $H^{***}$.

In the case the gluing is though a solid generalized cut, the dessin $D$ of a morsification is elementary equivalent to the gluing $D'-I_{2}^{**}$, where $I_{2}^{**}$ is the non-hyperbolic cubic corresponding to a line and a non-singular conic in $\RPP$ intersecting in two different real points, after a strict transform with respect to the blow up of a point in $\RPP$ which does not lie on any component.
Then the wiring diagram of $D$ would have either $1212$ or $2121$ corresponding to one of the nodal {\tvs} of $D'$ and the nodal {\tvs} of the gluing with $I_{2}^{**}$.

With this property in mind, we define the following operation on sequences.
Given a sequence $(a_i)_{i=1}^{3n}$, we put
\[
\eta_{\circ}((a_i)_{i=1}^{3n})=\left\{
\begin{array}{ll}
(a_1,a_2,\dots,a_n,1,2,1) & \text{ if $a_n=1$,}\\
(a_1,a_2,\dots,a_n,2,1,2) & \text{ if $a_n=2$,}\\
\end{array}\right.
\]

\[
\eta_{\times}((a_i)_{i=1}^{3n})=\left\{
\begin{array}{ll}
(a_1,a_2,\dots,a_n,2,1,2) & \text{ if $a_n=1$,}\\
(a_1,a_2,\dots,a_n,1,2,1) & \text{ if $a_n=2$.}\\
\end{array}\right.
\]

These function are not invariant under the action of $\sigma_{3n}$.

\begin{thm} \label{thm:wir}
Every wiring diagram $\omega$ of a morsification can be obtained as the projective equivalent class of a sequence of operations $\eta_{\circ}$, $\eta_{\times}$ and shifts of the wiring diagram $(1,2,1)$, i.e., 
there exist $\eta_i\in\{\eta_{\circ}, \eta_{\times}\}$, $k_i\in\N$, $i=1,2,\dots {n-1}$ such that $\omega$ is the projective equivalence class of the sequence
$$\eta_{n-1}\circ\sigma_{3n-3}^{k_{n-1}}\circ\dots\circ\eta_{2}\circ\sigma_{6}^{k_2}\circ\eta_1\circ\sigma_{3}^{k_1}(1,2,1).$$
\end{thm}

\begin{proof}
Let $C$ be a trigonal morsification such that $\omega$ is its wiring diagram.
Let $D=\Dssn(C)$ be its dessin.
By Proposition \ref{prop:decomptoi}, the dessin $D$ is elementary equivalent to a dessin having a decomposition $D_{3n-3}-D_{3}$, where $D_{3n-3}$ and~$D_3$ are dessins of degree ${3n-3}$ and $3$, respectively.
Fix a white vertex $v\in\Ver(D)$ and an orientation on the real part of $D$ such that $\omega$ is written as the wiring diagram of~$D$ with respect to $v$.

Let $k_{n-1}$ the number of nodal {\tvs} in the arc $(v,u)$, where $u$ is the vertex on the generalized cut passing from $D_3$ to $D_{3n-3}$ with respect to the aforementioned order of the real part of the dessin $D$.

Let $\eta_{n-1}$ be the operation $\eta_{\circ}$ if $D_3$ is the cubic dessin $H^{***}$ or the operation $\eta_{\times}$ if $D_3$ is the cubic dessin $I_{2}^{**}$.
In the latter case, we contract the solid segment of $D_{3n-3}$ to obtain the dessin $D_{3n-3}'$ of a morsification of degree $3n-3$. 
We iterate this decomposition $n-1$ times.
In the eventual case where the last cubic dessin has a resulting wiring diagram $(2,1,2)$, we add one to $k_1$ in order to obtain the desire decomposition.
\end{proof}

\begin{coro}
If $C$ is a morsification of degree $3n$,
then, all the entries on its cyclic length vector are less or equal to $n$.

Moreover, if its cyclic length vector has $n$ as an entry, then, the wiring diagram of $C$ is the only projective equivalence class with cyclic length vector 
\[
(1,n,1,\underbrace{2,2,\dots,2}_{{n-1}\text{ times.}})
\]
\end{coro}

\begin{proof}
We proceed by induction on $n$.
For $n=1$, there is a unique deformation class of morsifications in $\Sigma_1$, corresponding to the strict transform of three lines in $\RPP$ in general position, after the blow-up of a point which does not lie in any line.

Let us assume that $n>1$ and the statement holds for $n-1$.
Let $\omega$ be the wiring diagram of $C$.
By Theorem \ref{thm:wir}, we have that $\omega$ can be obtained as the projective equivalence class of $\eta_{n-1}(a_i)_{i=1}^{3n-3}$, with $(a_i)_{i=1}^{3n-3}=\sigma_{3n-3}^{k_{n-1}}\circ\dots\circ\eta_{2}\circ\sigma_{6}^{k_2}\circ\eta_1\circ\sigma_{3}^{k_1}(1,2,1).$
Denote by $o$ the cyclic length vector of $\omega$ and by $o_{n-1}$ the cyclic length vector of the projective equivalence class of $(a_i)_{i=1}^{3n-3}$.

In the case when $\eta_{n-1}=\eta_{\circ}$, we have that $o(a_{3n-3})$, the entry in $o$ corresponding to $a_{3n-3}$, increases at most by $1$ since $\eta_{\circ}$ concatenates to the sequence exactly one more element that equals $a_{3n-3}$.
Also, the entry $o(a_1)\leq o_{n-1}(a_1)+1$. Equality holds when $a_1=a_{3n-3}$ and $n-1$ is odd or when $a_1\neq a_{3n-3}$ and $n-1$ is even.
Every other entry in $o$ either comes from an entry in $o_{n-1}$ unaffected by $\eta_{n-1}$ or is the entry $1$ corresponding to the middle number in the added elements.
Therefore, since every entry in $o_{n-1}$ is at most $n-1$, then every entry in $o$ is at most $n$.

Moreover, if $o(a_{3n-3})$ or $o(a_1)$ equals $n$, then $o_{n-1}(a_{3n-3})$ or $o_{n-1}(a_1)$ equals ${n-1}$.
By the induction hypothesis, we have that $a_{3n-3}, a_1$ are distinct for $n-1$ even, or equal for $n-1$ odd and that 
either $o_{n-1}(a_{3n-3})=n-1$, $o_{n-1}(a_1)=1$,  
or $o_{n-1}(a_{3n-3})=1$, $o_{n-1}(a_1)=n-1$, respectively.
Therefore, in both cases $o(a_{3n-3})=o_{n-1}(a_{3n-3})+1$ and $o(a_1)=o_{n-1}(a_1)+1$.
Thus, the cyclic length vector $o$ equals $(1,n,1,2,2,\dots,2)$

In the case when $\eta_{n-1}=\eta_{\times}$, we have that $o(a_{3n-3})$ cannot increase since $\eta_{\times}$ concatenates to the sequence an element that differs from $a_{3n-3}$.
The entry ${o(a_1)\leq o_{n-1}(a_1)+1}$. Equality holds when $a_1=a_{3n-3}$ and $n-1$ is even or when $a_1\neq a_{3n-3}$ and $n-1$ is odd.
Every other entry in $o$ either comes from an entry in $o_{n-1}$ unaffected by $\eta_{n-1}$ or is the entry $1$ corresponding to the middle number in the added elements.
Therefore, since every entry in $o_{n-1}$ is at most $n-1$, then every entry in $o$ is at most $n$.

We cannot have that $o(a_1)=n$. Indeed, if $o(a_1)=n$, then $o_{n-1}(a_1)=n-1$ and $a_1=a_{3n-3}$ if $n-1$ is even or $a_1\neq a_{3n-3}$ if $n-1$ is odd.
On the other hand, since $o(a_1)=n$ and only the last element in $\eta_{\times}(a_i)_{i=1}^{3n-3}$ add to it, then $a_1=a_2=\dots=a_{n-1}$.
Then, the fact that $a_1=a_{3n-3}$ if $n-1$ is even or $a_1\neq a_{3n-3}$ if $n-1$ is odd, implies that the entry $o_{n-1}(a_1)\geq n$, which is a contradiction.
\end{proof}

It was remarked that in degree~$3$ the wiring diagram $(1,1,1,2,2,2,1,1,1)$ does not exist. Indeed, its cyclic length vector being $(3,3,3)$ does not correspond to a morsification.

%% file: mor2c.tex
\begin{figure}[t]
\begin{subfigure}[b]{0.4\textwidth}
\includegraphics[width=\textwidth]{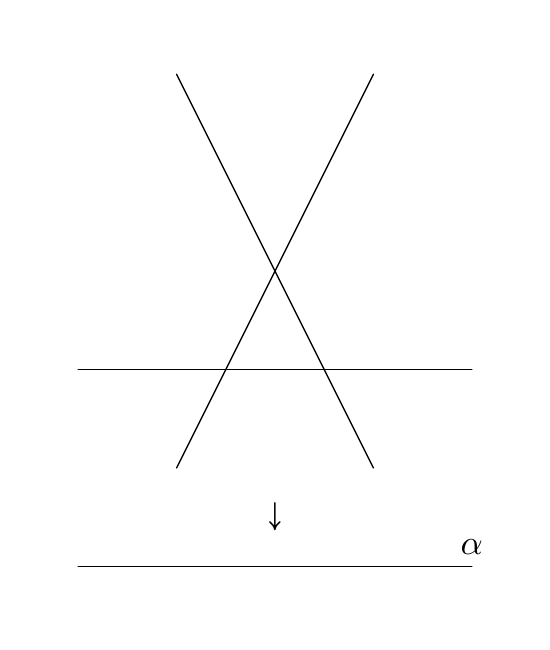}
\end{subfigure}
\begin{subfigure}[b]{0.4\textwidth}
\includegraphics[width=\textwidth]{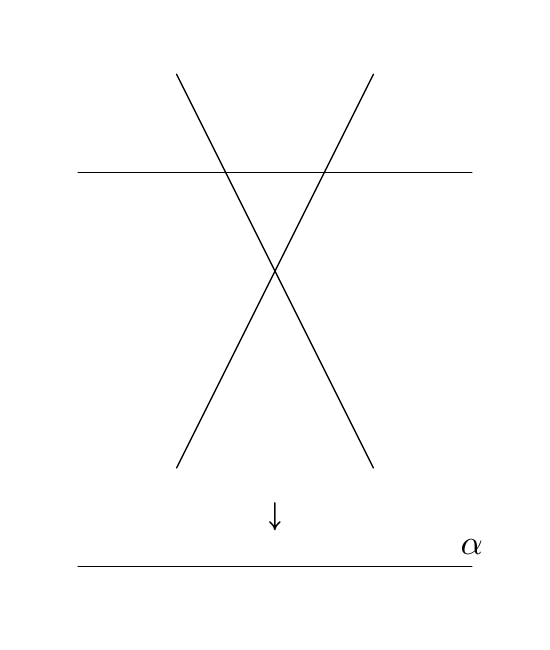}
\end{subfigure}
\caption{}
\label{fig:Reide}
\end{figure}

\begin{df}
Let $M\subset\Sigma_n$ be a trigonal morsification and let $X\subset\R\Sigma_n\lra\RP$ be a topological subspace.
If there exists an arc $\alpha\subset\RP$ such that $(\R\Sigma_n |_{\RP\setminus\alpha},X)$ is homeomorphic to $(\R\Sigma_n |_{\RP\setminus\alpha},M)$, and $(\R\Sigma_n |_{\alpha},X)$, $(\R\Sigma_n |_{\alpha},M)$ are homeomorphic to the curves in Figure~\ref{fig:Reide}, we say that $X$ \emph{is obtained from $M$ by a Reidemeister move}.

A trigonal morsification $M'\subset\Sigma_n$ is said \emph{to be obtained from $M$ by a Reidemeister move} if $\R M\subset\R\Sigma_n$ is.

We say that two trigonal morsifications $M, M'\subset\Sigma_n$  are \emph{Reidemeister equivalent} if there exists a chain $M_0=M, M_1,\dots, M_k=M'$ of trigonal morsifications such that $M_{i}$ is obtained from $M_{i-1}$ by a Reidemeister move, for $i=1,2, \dots,k$.
\end{df}

\begin{figure}[t]
\includegraphics[width=0.6\textwidth]{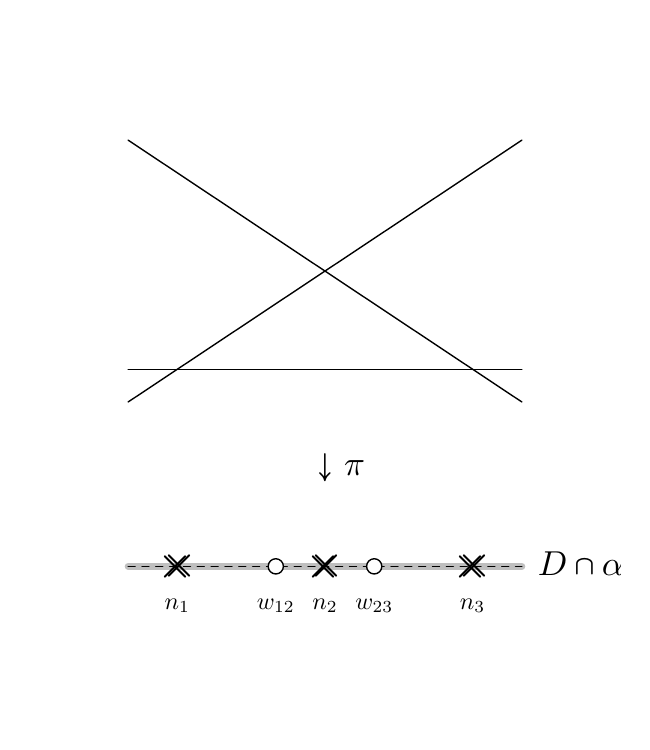}
\caption{Decorations of the segment of the dessin corresponding to the arc of a Reidemeister move.} \label{fig:ReideD}
\end{figure}

\begin{coro}
Let $M\subset\Sigma_n$ be a trigonal morsification.
If $X\subset\R\Sigma_n$ is obtainbed from $M$ by a Reidemeister move, then, there exists a trigonal morsification $M'\subset\Sigma_n$ such that $(\R\Sigma_n,X)\cong(\R\Sigma_n,\R M')$.
\end{coro}

\begin{figure}[b]
\begin{subfigure}[b]{0.45\textwidth}
\includegraphics[width=\textwidth]{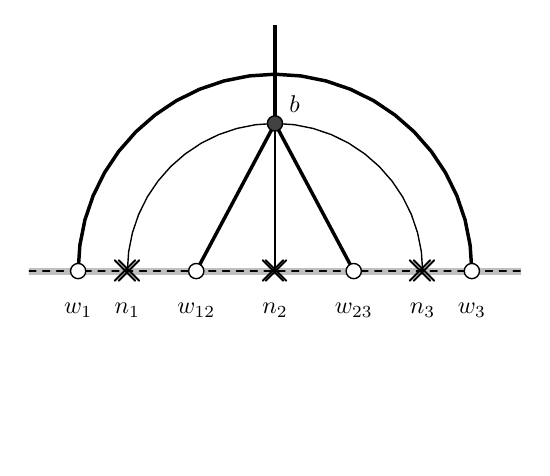}
\caption{}
\label{fig:CoroA}
\end{subfigure}
\begin{subfigure}[b]{0.45\textwidth}
\includegraphics[width=\textwidth]{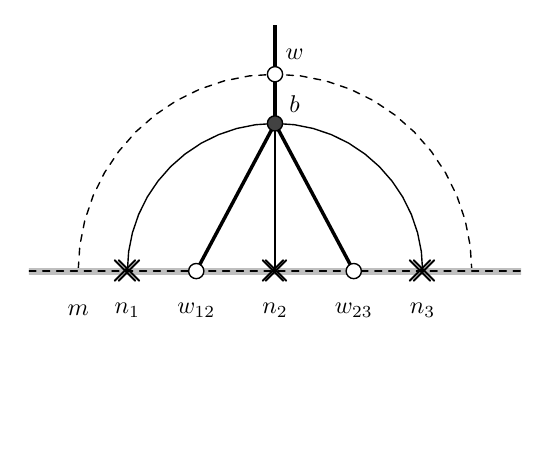}
\caption{}
\label{fig:CoroB}
\end{subfigure}
\caption{}\label{fig:Coro}
\end{figure}

\begin{proof}
Let us consider $n\geq 2$ since the statement holds trivially for $n=1$.
Let $D\subset\CP/(z\lmt\bar{z})$ be the dessins associated to $M$.
Let us consider the segment $S_{\alpha}\subset D$ corresponding to the arc $\alpha\subset\RP$ such that $(\R\Sigma_n |_{\alpha},X)$ and $(\R\Sigma_n |_{\alpha},M)$ are homeomorphic to the curves in Figure~\ref{fig:Reide}.

We can assume $D$ is bridge-free (see \cite[Lemma 2.3]{deg5}).
Let us denote by $n_1, n_2, n_3$ the {\tvs} in $S_\alpha$ corresponding to the nodal vertices.
Up to elementary moves of type $\circ$-in, we can assume there is exactly one real white vertex $w_{12}$ neighbouring $n_1$ and $n_2$ and exactly one real white vertex $w_{23}$ neighbouring $n_2$ and $n_3$ 
(see Figure~\ref{fig:ReideD}).
Up to monochrome modification, the vertices $n_1, w_{12}, n_2, w_{23}$ and $n_3$ are adjacent to an inner black vertex $b$.

If $n_1$ and $n_3$ have neighboring white vertices $w_1\neq w_{12}$ and $w_3\neq w_{23}$, the creation of an inner bold monochrome vertex adjacent to $w_1$, $b$ and $w_3$ bring us to the configuration shown in Figure~\ref{fig:CoroA}.

Otherwise, if there is a monochrome vertex $m$ adjacent to $n_1$ or $n_3$, up to an elementary move of type $\circ$-in, the vertex $m$ is connected to an inner white vertex~$w$.
Up to a monochrome modification, the vertex $w$ is adjacent to $b$.
Then, up to the creation of a dotted bridge or a monochrome modification, the dessin $D$ has a subgraph shown in Figure~\ref{fig:CoroB}.

In each case, replacing this subgraph with the alternative one in Figure~\ref{fig:Coro} produces a dessin $D'$ whose correponding trigonal curve $M'\subset\Sigma_n$ is a morsification such that $(\R\Sigma_n,X)\cong(\R\Sigma_n,\R M')$ since outside the segment $S_\alpha$ there are no changes and within the segment $S_\alpha$ the parity of the number of white vertices in the dotted segments neighboring $n_1$ and $n_3$ changed.
\end{proof}

\begin{coro}
If $M, M'\subset\Sigma_n$ are trignal moorsatificions, then, they are Reidemeister equivalent.
\end{coro}

\begin{proof}
We prove that every trigonal morsification $M\subset\Sigma_n$ is Reidemeister equivalent to the trigonal morsification $M^*$ whose wiring diagram is
$$\omega_n^{*}:=\eta_{\times}^{n-1}(1,2,1).$$

Let us remark that $\omega_n^*$ is invariant under the action of $\sigma_{3n}$, for all $n$, and the fact that a Reidemeister move on a morsification change a block of entries of its wiring diagram of the form $(1,2,1)$ for a block of the form $(2,1,2)$ and vice versa.
Therefore, for every $i$, we have that the morsifications corresponding to 
$\eta_i\circ\sigma_i^{k_i}\omega_i^*$ and $\omega_i^*$ are Reidemeister equivalent.

The statement follows by applying these facts to the decomposition 
$\eta_{n-1}\circ\sigma_{3n-3}^{k_{n-1}}\circ\dots\circ\eta_{2}\circ\sigma_{6}^{k_2}\circ\eta_1\circ\sigma_{3}^{k_1}(1,2,1)$
of the wiring diagram $\omega$ of $M$ given by Theorem~\ref{thm:wir}.\end{proof}

%% file: appendix2.tex
\newcommand{\Int}{{\operatorname{Int}}}\newcommand{\MCTI}{{\operatorname{MCTI}}}
\newcommand{\AG}{{\mathrm A}\Gamma}\newcommand{\MCTE}{{\operatorname{MCTE}}}
\newcommand{\Ima}{{\operatorname{Im}}}
\newcommand{\Rea}{{\operatorname{Re}}}
\newcommand{\BH}{{\operatorname{bh}}}
\newcommand{\Span}{{\operatorname{Span}}}
\newcommand{\ini}{{\operatorname{ini}}}
\newcommand{\fm}{{\mathfrak{m}}}
\newcommand{\Tor}{{\operatorname{Tor}}}
\newcommand{\conv}{{\operatorname{Conv}}}
\newcommand{\Conj}{{\operatorname{conj}}}
\newcommand{\Ev}{{\operatorname{Ev}}}
\newcommand{\pr}{{\operatorname{pr}}}
\newcommand{\val}{{\operatorname{val}}}
\newcommand{\eps}{{\varepsilon}}
\newcommand{\DD}{\boldsymbol{D}}
\newcommand{\idim}{{\operatorname{idim}}}
\newcommand{\DP}{{\operatorname{DP}}}
\newcommand{\Arc}{{\operatorname{Arc}}}
\newcommand{\Graph}{{\operatorname{Graph}}}
\newcommand{\sm}{{\mathrm{sm}}}
\newcommand{\bw}{{\boldsymbol{w}}}
\newcommand{\D}{{\mathbb{D}}}
\newcommand{\bz}{{\boldsymbol{z}}}
\newcommand{\bp}{{\boldsymbol{p}}}
\newcommand{\const}{{\operatorname{const}}}
\newcommand{\Tors}{{\operatorname{Tors}}}
\newcommand{\Spec}{{\operatorname{Spec}}}
\newcommand{\bA}{{\boldsymbol{A}}}\newcommand{\bF}{{\boldsymbol{F}}}
\newcommand{\bB}{{\boldsymbol{B}}}
\newcommand{\bn}{{\boldsymbol{n}}}
\newcommand{\bs}{{\boldsymbol{s}}}
\newcommand{\PP}{{\mathbb{P}}}
\newcommand{\codim}{{\operatorname{codim}}}
\newcommand{\Ann}{{\operatorname{Ann}}}
\newcommand{\ord}{{\operatorname{ord}}}
\newcommand{\mt}{{\operatorname{mt}}}\newcommand{\Br}{{\operatorname{Br}}}
\newcommand{\Prec}{{\operatorname{Prec}}}\newcommand{\ReBr}{{\operatorname{ReBr}}}
\newcommand{\ImBr}{{\operatorname{ImBr}}}
\newcommand{\proofend}{\hfill$\blacksquare$\bigskip}
\newcommand{\pperp}{{\perp{\hskip-0.4cm}\perp}}
\newcommand{\bk}{{\boldsymbol{k}}}
\newcommand{\bl}{{\boldsymbol{l}}}
\newcommand{\cP}{{\mathcal P}}

\begin{appendix}

\section{E. Shustin. Polynomiality of morsifications of trigonal singularities}
\label{app}
%

\subsection{Main result}

By a {\bf singularity} we always mean a germ $(C,z)\subset\C^2$ of a plane reduced analytic curve at its singular point $z$.
Irreducible components of the germ
$(C,z)$ are called {\bf branches of $(C,z)$}. Let $f(x,y)=0$ be an (analytic) equation of $(C,z)$, where
$f$ is defined in the closed ball $B(z,\eps)\subset\C^2$
of radius $\eps>0$ centered at
$z$. The ball $B(z,\eps)$ is called the {\bf Milnor ball} of $(C,z)$ (and is denoted in the sequel $B_{C,z}$) if
$z$ is the only singular point of $C$ in $B(z,\eps)$, and $\partial B(z,\eta)$ intersects $C$ transversally
for all $0<\eta\le\eps$. A {\bf nodal deformation} of a singularity $(C,z)$ is a
family of analytic curves $C_t=\{f_t(x,y)=0\}$, where $f_t(x,y)$ is analytic in $x,y,t$ for $(x,y)\in B(C,z)$ and
$t$ varying in an open disc $\D_\zeta\subset\C$ of some radius $\zeta>0$ centered at zero, and where $C_0=C$, $C_t$ is smooth along $\partial B_{C,z}$, intersects $\partial
B_{C,z}$ trasversally for all $t\in\D_\zeta$,  the curve $C_t$ has only ordinary nodes in $B_{C,z}$ for any $t\ne0$,
and the number of nodes does not depend on $t$.
The maximal number of nodes in a nodal deformation of $(C,z)$ in $B$ equals $\delta(C,z)$, the $\delta$-invariant
(see, for instance, \cite[\S10]{M}).

Let $(C,z)$ be a totally real singularity, i.e.,
invariant with respect to the complex conjugation, with $z\in C$ a real singular point and
all branches real. Let $C_t=\{f_t(x,y)=0\}$, $t\in\D_\zeta$, be an equivariant\footnote{Here
and further on, {\it equivariant} means commuting with the complex conjugation.} nodal deformation
of a real singularity $(C,z)$. Its restriction to $t\in[0,\zeta)$ is called
a {\bf real nodal deformation}.
A real nodal deformation is called a {\bf real morsification} of $(C,z)$ if
any $C_t$, $t>0$, has $\delta(C,z)$ hyperbolic real nodes (i.e., intersection points of two smooth transverse real branches).

Denote by $T_\R(3,3n)$ the class of real plane curve singularities having three
real smooth branches intersecting each other with multiplicity $n$.


\begin{thm}\label{tb}
Let $(C,0)\subset B(0,\eps)\subset\C^2$ be a germ of a real plane curve singularity of type $T_\R(3,3n)$, given by
an analytic equation
$$F(x,y)=f_{3n}+\sum_{3i+nj>3n}a_{ij}x^iy^j=0,\quad f_{3n}=y^3+a_{n,2}x^ny^2+a_{2n,1}x^{2n}y+a_{3n,0}x^{3n}\ ,$$
where $f_{3n}$ is squarefree, and let
$F_t(x,y)$, $|t|<\eta$, be a morsification of the given singularity
$(C,0)$ (respectively having $3n$ hyperbolic nodes in $B^\R(0,\eps)$ for each $0<t<\eta$). Then there exists a polynomial
\begin{equation}G(x,y)=f_{3n}+\sum_{3i+nj<3n}a_{ij}x^iy^j\label{eb}\end{equation} such that $\{G=0\}\cap B(0,\eps)$ is
equivariantly isotopic
to any  curve $\{F_t=0\}\cap B(0,\eps)$, $|t|<\eta$, relative to $\partial B(0,\eps)$.
\end{thm}

In the proof, we explore the idea originated in the proof of \cite[Corollary]{Sh1}, later used in
\cite[Proof of Theorem 3]{KS}, and elaborated further in \cite[Proof of Theorem 5]{OS}. It, however, appears to be
rather more involved technically than \cite[Proof of Theorem 5]{OS}, and we present it in several parts.

\begin{rmk}\label{rb}
The statement of Theorem \ref{tb} yields the following fact: if $(C,z)$ is a real quasihomogeneous singularity of
type $T_\R(3,3n)$, and $(C',z)$ is a real quasihomogeneous singularity of type $T_\R(3,3n)$ obtained from $(C,z)$ by a sufficiently small
equisingular deformation, then the set of isotopy types of morsifications of $(C,z)$ is
naturally included into the set of isotopy types of morsifications of $(C',z)$. Indeed, if, in the notation
of Theorem \ref{tb}, a morsification
of $(C,z)=\{f_{3n}=0\}$ is isotopic to $\{G=0\}\cap B(z,\eps)$, where
$$G(x,y)=\prod_{i=1}^3\left(y+\sum_{j=0}^n\alpha_{ij}x^j\right)\ ,$$ then a small variation of the factors of $G(x,y)$ preserves all
real nodes of the curve $\{G=0\}$, and hence its equivariant isotopy type.
\end{rmk}

\subsection{Proof of Theorem \ref{tb}}

\subsubsection{Part I}

We start with reducing the problem to the independence of simultaneous morsifications
of singularities of a real curve given by a polynomial with the the Newton polygon under
the line $\{i+nj=3n\}$.

\smallskip

{\bf(1)} Observe that there exits $\alpha$ such that after the coordinate change
$(x,y)\mapsto(x,y+\alpha x^n)$, the coefficients of $f_{3n}$ satisfy the relation
\begin{equation}9a_{3n,0}-a_{n,2}a_{2n,1}\ne0\ .\label{eb1}\end{equation}
Indeed, the above coordinate change turns the left-hand side of (\ref{eb1}) into
$$2\alpha(3a_{2n,1}-a_{n,2}^2)+9a_{3n,0}-a_{2n,1}a_{n,2}\ ,$$ which vanishes identically only if $f_{3n}$
is the cube of a binomial against our assumptions. Since the coordinate change we used does not affect
the statement of proposition, we can assume that (\ref{eb1}) holds.

\smallskip

{\bf(2)} For a given
\begin{equation}F_t(x,y)=\sum_{i,j\ge0}b_{ij}(t),\quad b_{ij}(t)=\begin{cases}a_{ij},\quad & i+nj\ge3n,\\
0,\quad & i+nj<3n,\end{cases}\label{eb2}\end{equation} there exist analytic functions $\alpha(t),\beta_0(t),...,\beta_{n-1}(t)$
vanishing at zero such that $F_t(x+\alpha(t),y+\beta_0(t)+\beta_1(t)x+...+\beta_{n-1}(t)x^{n-1})$ does not contain
the monomials $y^2,...,x^{n-1}y^2$, and $x^{3n-1}$. Indeed, the required property reduces to a system of
equations on $\alpha,\beta_1,...,\beta_{n-1}$ with the Jacobian at $t=0$ equal to
$3^{n-1}(9a_{3n,0}-a_{n,2}a_{2n,1})$, non-vanishing according to (\ref{eb1}). Thus, without loss of generality,
we can assume that
$$b_{0,2}(t)=...=b_{n-1,2}(t)=b_{3n-1,0}(t)\equiv0\ .$$ Since $F_t=0$ has no singularity of type
$T_\R(3,3n)$ as $t\ne0$, we get $\sum_{i+nj<3n}|b_{ij}(t)|>0$ for $t\ne0$. Hence, for a fixed
$c_0>0$, there exists a positive function $\tau(t)$ on the interval $(0,\eta)$ such that
\begin{equation}\sum_{i+nj<3n}|b_{ij}(t)|\tau(t)^{3n-i-nj}=c_0,\quad t>0\ .\label{eb3}\end{equation} In view of (\ref{eb2}), $\lim_{t\to0}\tau(t)=\infty$.
Thus, in the family
$$\tau(t)^{3n}F_t(x\tau(t)^{-1},y\tau(t)^{-n})=
\sum_{i+nj>3n}b_{ij}(t)\tau^{3n-i-nj}x^iy^j$$
\begin{equation}+\sum_{i+nj=3n}b_{ij}(t)x^iy^j
+\sum_{i+nj<3n}b_{ij}(t)t^{3n-i-nj}x^iy^j\ ,\label{eb4}\end{equation}
the first sum converges to zero, while the second one converges to $f_{3n}(x,y)$ as $t\to0$.
Moreover, by (\ref{eb3}), there exists a sequence $t_m\to0$, $m=1,2,...$, such that
the third sum in (\ref{eb4}) converges as well so that
$$\lim_{m\to\infty}\tau(t_m)^{3n}F_t(x\tau(t_m)^{-1},y\tau(t_m)^{-n})=f_{3n}(x,y)+\sum_{i+nj<3n}a_{ij}x^iy^j=:R(x,y)\ ,$$
where
$$\sum_{i+nj<3n}|a_{ij}|=c_0>0\quad\text{and}\quad a_{0,2}=...=a_{n-1,2}=a_{3n-1,0}=0\ .$$\emph{}
Suppose that $c_0>0$ is chosen so that the curve $\{R(x,y)=0\}\subset\C^2$ is smooth outside $B(0,\eps)$ and
intersects with all spheres $\partial B(0,\zeta)$, $\zeta\ge\eps$, transversally.
Thus, by construction we have:
\begin{enumerate}\item[(R1)] $\{R(x,y)=0\}$ has no singularity of type $T_\R(3,3n)$ in $B(0,\eps)$,
and each its singular point is a center of two or three smooth real branches; furthermore, since
the singularities of $\{R=0\}$ admit a deformation into $3n$ nodes in total, $R(x,y)$
splits into the product of $(y-Q_1(x))(y-Q_2(x))(y-Q_3(x))$, where $\deg Q_i=n$, $i=1,2,3$;
\item[(R2)] there exists an equivariant analytic deformation of $\{R(x,y)=0\}$, simultaneously realizing
morsifications of all singular points of $\{R=0\}$ so that the resulting curve in $B(0,\eps)$
is equivariantly isotopic to
$\{F_t(x,y)=0\}\cap B(0,\eps)$, $t>0$, relative to $\partial B(0,\eps)$.
\end{enumerate}
Denote by $\cP(k)\subset\C[x,y]$ the linear subspace spanned by the monomials $x^iy^j$ with
$i+nj=k$, and put
$$\cP_{(k)}=\bigoplus_{i\le k}\cP(i),\quad \cP^{(k)}=\bigoplus_{i\ge k}
\cP(i)\ .$$ So, we will complete the proof of Theorem \ref{tb} when showing that
the germ at $R$ of the affine space $R+\cP_{(3n-1)}$ induces one-parameter deformations
simultaneously realizing morsifications of prescribed isotopy types for all singular points
of the curve $\{R=0\}\cap B(0,\eps)$.
We shall prove this fact by induction on $n$. For $n=1$, we have a $D_4$ singularity, for which the
statement is evident. The induction step splits into several cases treated below.

Before proceeding further on, we make one more coordinate change $(x,y)\mapsto(x,y+\alpha x^n)$, which
annihilates the coefficient of $x^ny^2$, i.e., leads to
\begin{equation}\sum_{i+nj<3n}|a_{ij}|=c'_0>0\quad\text{and}\quad a_{0,2}=...=a_{n-1,2}=a_{n,2}=0\ ,
\label{eb5a}\end{equation} while, of course, preserving the properties (R1), (R2).

\subsubsection{Part II}

Suppose that
$R(x,y)=yS(x,y)$. In view of (\ref{eb5a}), $S(x,y)=y^2-Q(x)^2$, where
a polynomial $Q(x)$ of degree $n$ has at least two roots and all the roots are real.
That is, the curve $C(R):=\{R=0\}$ has singularities of types $T_\R(3,3n_i)$, $1\le i\le k$ ($k\ge2$), where
$n_1+...+n_k=n$. For each singular point $z_i$ of $C(R)$, $i=1,...,k$, denote by
$I_i\subset{\mathcal O}_{\C^2,z}$ the ideal defined by vanishing of the coefficients of the monomials
lying strictly below the line $i+n_ij=3n_i$. Then the statement of the induction step can be reformulated
as the surjectivity of the natural map $f_{3n}+\cP_{3n-1}\to\bigoplus_i{\mathcal O}_{\C^2,z}/I_i$.

To prove this surjectivity, we consider $C(R)$ as a
curve in the linear system $|3D|$ on the Hirzebruch surface ${\mathcal F}_n$,
where $D$ is the divisor class of the section disjoint from the exceptional curve. Introduce the zero-dimensional subscheme $Z\subset{\mathcal F}_n$,
$Z=Z_0\cup Z_1\cup...\cup Z_k$, where $Z_0$ consists of the three reduced points of intersection of $C(R)$ with
the fiber $x=\infty$, and $Z_i$ is defined at the point $z_i$ by the ideal $I_i$, $i=1,...,k$. Then the required statement
admits the following cohomological reformulation:
\begin{equation}H^1({\mathcal F}_n,{\mathcal J}_{Z/{\mathcal F}_n}(3D))=0\ ,\label{eb8}\end{equation}
where ${\mathcal J}_Z$ is the ideal sheaf of the subscheme $Z$. The curve $C(R)$ splits into three components
$C_1,C_2,C_3\in|D|$. We apply ``le methode d'Horace" (see, \cite{Hir}). Consider the three exact sequences
$$0\to{\mathcal J}_{Z:C_1/{\mathcal F}_n}(2D)\overset{C_1}{\to}{\mathcal J}_{Z/{\mathcal F}_n}(3D)
\to{\mathcal J}_{Z\cap C_1/C_1}(3D)\to0\ ,$$
$$0\to{\mathcal J}_{Z:(C_1C_2)/{\mathcal F}_n}(D)\overset{C_2}{\to}{\mathcal J}_{Z:C_1/{\mathcal F}_n}(2D)
\to{\mathcal J}_{(Z:C_1)\cap C_2/C_2}(2D)\to0\ ,$$
$$0\to{\mathcal O}_{{\mathcal F}_n}\overset{C_3}{\to}{\mathcal J}_{Z:(C_1C_2)/{\mathcal F}_n}(D)
\to{\mathcal J}_{(Z:(C_1C_2))/C_3}(D)\to0\ ,$$ which yield the exact cohomology sequences
$$H^1({\mathcal F}_n,{\mathcal J}_{Z:C_1/{\mathcal F}_n}(2D))\to
H^1({\mathcal F}_n,{\mathcal J}_{Z/{\mathcal F}_n}(3D))\to
H^1(C_1,{\mathcal J}_{Z\cap C_1/C_1}(3D))\ ,$$
$$H^1({\mathcal F}_n,{\mathcal J}_{Z:(C_1C_2)/{\mathcal F}_n}(D))\to
H^1({\mathcal F}_n,{\mathcal J}_{Z:C_1/{\mathcal F}_n}(2D))\to
H^1(C_2,{\mathcal J}_{(Z:C_1)\cap C_2/C_2}(2D))\ ,$$
$$0=H^1({\mathcal O}_{{\mathcal F}_n}\to H^1({\mathcal F}_n,{\mathcal J}_{Z:(C_1C_2)/{\mathcal F}_n}(D))\to
H^1(C_3,{\mathcal J}_{(Z:(C_1C_2))/C_3}(D))\ ,$$
(here $Z:(C_1C_2)\subset C_3$) and hence (\ref{eb8}) can be derived from
$$H^1(C_1,{\mathcal J}_{Z\cap C_1/C_1}(3D))=
H^1(C_2,{\mathcal J}_{(Z:C_1)\cap C_2/C_2}(2D))=
H^1(C_3,{\mathcal J}_{(Z:(C_1C_2))/C_3}(D))=0\ .$$ To establish these $h^1$-vanishing relations, we
use 
the Riemann-Roch criterion, which for the smooth rational curves $C_1,C_2,C_3$ respectively reads
$$3D^2-\deg(Z\cap C_1)>D^2+DK_{{\mathcal F}_n}=-2,\quad 2D^2-\deg((Z:C_1)\cap C_2)>-2\ ,$$
$$D^2-\deg(Z:(C_1C_2))>-2\ ,$$ and all these conditions are fulfilled in view of $D^2=n$ and
$$\deg(Z\cap C_1)=3n+1,\quad \deg((Z:C_1)\cap C_2)=2n+1,\quad \deg(Z:(C_1C_2))=n+1\ .$$

\subsubsection{Part III}

Suppose that $R(x,y)$ is not divisible by $y$, and the curve $C(R)$ contains a singular point of multiplicity $3$.
It follows from the properties (R1), (R2), and (\ref{eb5a}) that this point $z=(x_0,y_0)$ must be of
type $T(3,3m)$ with $1\le m<n$, and that $R(x_0,y)=y^3$. So, after the shift $(x,y)\mapsto(x-x_0,y)$ followed by
an appropriate change $(x,y)\mapsto(x,y+\alpha_1x+...+\alpha_mx^m)$, we obtain a polynomial
$\widetilde R(x,y)$ with the Newton triangle $\conv\{(3m,0),(3n,0),(0,3\}$, whose truncations to the edges
$[(3m,0),(0,3)]$ and $[(3n,0),(0,3)]$ are squarefree. The polynomial $\widehat R(x,y):=x^{-3m}\widetilde R(x,x^my)$ has Newton trinagle $\conv\{(0,0),(3n-3m,0),(0,3)\}$ and defines a curve with the same singularities in $\C^2\setminus\{x=0\}$ as the
polynomial $\widetilde R(x,y)$. Thus, by the induction assumption, there exists a deformation of the coefficients of the monomials under
the segment $[(3n-3m,0),(0,3)]$, $\widehat R_t(x,y)$, $0\le t<\eta$, simultaneously realizing morsifications of prescribed isotopy types
for all the singularities of $\{\widehat R=0\}$. Hence the deformation
of $\widetilde R(x,y)$ given by $\widetilde R_t(x,y):=x^{2m}\widehat R_t(x,x^{-m}y)$, $0\le t<\eta$, simultaneously realizes
morsifications of prescribed isotopy types for all the singularities of $\{\widetilde R=0\}$, while keeping the singularity
of type $T_\R(3,3m)$  at the origin. By the induction assumption, Theorem \ref{tb} and Remark \ref{rb} apply to the singularity
at the origin of $\{\widetilde R_{t_0}=0\}$ (with some fixed $0<t_0<\eta$), which means that for any chosen
morsification of the singularity of $\{\widetilde R=0\}$ at the origin, there exists a real polynomial
$P(x,y)$ with Newton triangle $\conv\{(0,0),(3m,0),(0,3)\}$, the same truncation of
the edge $[(3m,0),(0,3)]$ as $\widetilde R_{t_0}(x,y)$, and such that the curve $\{P=0\}\cap B(0,\eps)$ is nodal
and equivariantly isotopic to the chosen morsification in $B(0,\eps)$ (relative to
$\partial B(0,\eps)$). By \cite[Theorems 3.1 and 4.1(1)]{ShT}, we can ``patchwork the polynomials $\widetilde R_{t_0}(x,y)$ and $P(x,y)$ and
obtain a family of real polynomials with Newton triangle $\conv\{(0,0),(3n,0),(0,3)\}$
simultaneously realizing morsifications of prescribed isotopy types for all the singularities of the original
curve $\{R=0\}$.

\subsubsection{Part IV}

So, we are left with the case of the curve $C(R)$ having only double singular points, i.e. points of intersection of
two local smooth real branches. 
 We shall show that, in this case, the linear system
$|3D|$ on the Hirzebruch surface ${\mathcal F}_n$ induces a joint versal deformation of all singular points of
the curve $C(R)\subset{\mathcal F}_n$.

\smallskip

{\bf(1)} It is well-known that a versal deformation of each point
$(x_0,y_0)\in\Sing(C(R))$ is generated by any basis of the linear space
$\C\{x-x_0,y-y_0\}/\langle R,R_x,R_y\rangle$,
where $\C\{*,*\}$ denotes the ring of locally convergent power series, and $\langle R,R_x,R_y\rangle$
is the Tjurina ideal (generated by $R,R_x,R_y$). Hence, it is enough to
prove the surjectivity of the following two projections:
$$\pr_1:\C[x,y]/\langle R,R_x,R_y\rangle\to\bigoplus_{(x_0,y_0)\in\Sing(C(R))}
\C\{x-x_0,y-y_0\}/\langle R,R_x,R_y\rangle\ ,$$
$$\pr_2:\cP_{3n}\to\C[x,y]/\langle R,R_x,R_y\rangle\ .$$
To prove the surjectivity of $\pr_1$,
introduce the zero-dimensional subscheme $Z$ of the Hirzebruch surface
${\mathcal F}_n$, concentrated
at $\Sing(C(R))$ and defined by the local Tjurina ideal at each point.
Then we have an exact sequence of sheaves
$$0\to{\mathcal J}_{Z/{\mathcal F}_n}(mD)\to{\mathcal O}_{{\mathcal F}_n}(mD)\to {\mathcal O}_Z\to 0\ ,$$
${\mathcal J}_{Z/{\mathcal F}_n}$ being the ideal sheaf of $Z$. Since for $m\gg0$,
$H^1({\mathcal J}_{Z/{\mathcal F}_n}(mD))=0$, we get
\begin{eqnarray}
     &H^0({\mathcal O}_{{\mathcal F}_n}(mD))/H^0({\mathcal J}_{Z/{\mathcal F}_n}(mD))
     \simeq H^0({\mathcal O}_Z)\qquad\qquad\nonumber\\
     &\qquad\qquad=\bigoplus_{(x_0,y_0)\in\Sing(\{R=0\})}
     \C\{x-x_0,y-y_0\}/\langle R,R_x,R_y\rangle,\nonumber
\end{eqnarray}
which, in fact, yields that $\pr_1$ is an isomorphism.

\smallskip

{\bf(2)} Now we prove the surjectivity of $\pr_2$.

Since $R(x,0)\not\equiv0$ and $C(R)$ is the union of three smooth real sections, there exist
$x_0,x_\infty\in\R$ such that $R(x_0,0)=0$ and $R(x_\infty,0)\ne0$. Performing an equivariant automorphism of the base of
${\mathcal F}_n$, we can set $x_0=0$, $x_\infty=\infty$, while keeping the property (\ref{eb5a}). Note, that then
$f_{3n}(x,y)=y^3+a_{2n,1}x^{2n}y+a_{3n,0}x^{3n}$ with $a_{2n,1}a_{3n,0}\ne0$, since $F_{3n}$ splits into three real factors.
Write
\begin{equation}R(x,y)=f_{3n}+\sum_{i=1}^sf_{k_i},\quad f_{k_i}\in\cP(k_i)\setminus\{0\},\quad
3n>k_1>...>k_s\ge0,\ s\ge1\ .\label{eb10}\end{equation}
Observe that $k_s>0$ and $f_{n,y}\ne0$, that is,
\begin{equation}a_{0,0}=R(0,0)=0\quad\text{and}\quad a_{0,1}=R_y(0,0)\ne0\ .
\label{eb20}\end{equation}
The former relation - by construction, the latter relation - due to (\ref{eb5a}) and the absence of triple singular points.

Since $f_{3n}$ is square-free, the
derivatives $f_{3n,x}\in\cP(3n-1)$ and $f_{3n,y}\in\cP(2n)$ are coprime. Hence, any polynomial
in $\cP^{(4n-1)}$ is congruent to a polynomial in $\cP_{(4n-2)}$ modulo $\langle R_x,R_y\rangle$.

Notice that $$f_{3n,x}=2na_{2n,1}x^{2n-1}y+3na_{3n,0}x^{3n-1}=nx^{2n-1}g,\quad
g=2a_{2n,1}y+3a_{3n,0}x^n\ .$$

To show that
\begin{equation}\cP{(m)}\equiv\cP_{(m-1)}\mod\langle R,R_x,R_y\rangle\quad\text{for all}\quad 3n\le m\le 4n-2\ ,
\label{eb14}\end{equation}
which is equivalent to the surjectivity of $\pr_2$, we construct a finite sequence of polynomials $R_i
\in\langle R,R_x,R_y\rangle$, $i\ge1$, such that the final polynomial $R_p(x,y)$ has the leading $(1,n)$-homogeneous form
\begin{equation}g_k=ax^{k-n}y+bx^k\not\in\langle g\rangle,\quad k<3n-1\ .\label{eb13}\end{equation}
Since (see (\ref{eb5a}))
$$f_{3n,y}=3y^2+a_{2n,1}x^{2n},\quad f_{3n,x}=2na_{2n,1}x^{2n-1}y+3na_{3n,0}x^{3n-1}\ ,$$
$$\cP(m)=\Span\{x^{m-3n}y^3,\ x^{m-2n}y^2,\ x^{m-n}x^{m-n}y,\ x^k\},\quad 3n\le m\le 4n-2\ .$$
relation (\ref{eb13}) will imply that
$$\cP(m)=\cP(m-2n)\cdot f_{3n,y}+\cP(m-3n+1)\cdot f_{3n,x}+\cP(m-k)\cdot g_k,\quad3n\le m\le 4n-2\ ,$$
and hence (\ref{eb14}).

\smallskip

We start with the required sequence $\{R_i(x,y)\}_{i\ge1}$. If $k_1=3n-1$ and $f_{k_1}\in\langle g\rangle$, then we
can perform a coordinate change $(x,y)\mapsto(x+\alpha,y)$, which annihilates the form $f_{3n-1}$, but does not affect
the properties (R1), (R2), and (\ref{eb5a}). Thus, we can assume that
\begin{equation}\text{either}\ f_{k_1}\not\in\langle g\rangle,\quad \text{or}\ k_1=3n-d,\ d\ge2,
\ f_{k_1}\in\langle g\rangle\ .\label{eb5b}\end{equation}

Euler's identities $3nf_{3n}=xf_{3n,x}+nyf_{3n,y}$ and $k_1f_{k_1}=xf_{k_1,x}+nyf_{k_1,y}$ yield 
\begin{equation}R_1(x,y):=3nR(x,y)-xR_x(x,y)-nR_y(x,y)=(3n-k_1)f_{k_1}+\sum_{i\ge2}(3n-k_i)f_{k_i}\ .\label{eb7}\end{equation}
The former option in (\ref{eb5b}) leads to
the situation of (\ref{eb13}) with $g_{k_1}=f_{k_1}$, and hence to
(\ref{eb14}) as required.

So, assume now that $f_{k_1}\in\langle g\rangle$, $k_1=3n-d$, $d\ge2$. Starting with the polynomials
$R_0:=R_x$ and $R_1$, we perform the following operation. Suppose that, for some $i\ge0$, we have a couple of
polynomials $R_i(x,y),R_{i+1}(x,y)$, containing only monomials of type $x^k$, $x^ly$, having the
highest $(1,n)$-quasihomegenelous forms $\lambda_ix^{m_i}g$ and $\lambda_{i+1}x^{m_{i+1}}g$, respectively,
where $\lambda_i\lambda_{i+1}\ne0$.
If, say, $m_i\ge m_{i+1}$, we define $R_{i+2}=\lambda_{i+1}R_i-\lambda_ix^{m_i-m_{i+1}}R_i$. If the leading form
of $R_{i+2}$ does not belong to $\langle g\rangle$, we are done (se the above argument). If $R_{i+2}\not\equiv0$, and the leading form of
$R_{i+2}$ belongs to $\langle g\rangle$, then we arrive to the pair $R_{i+1},R_{i+2}$ with the total
degree $m_{i+1}+m_{i+2}<m_i+m_{i+1}$. In the worst case, we end up with some $R_m(x,y)\equiv0$, which means that 
$$R_x(x,y)=P_1(x)Q(x,y),\quad R_1=P_2(x)Q(x,y)\ ,$$ where 
$$Q(x,y)=x^mg(x,y)+\sum_{\renewcommand{\arraystretch}{0.6}
\begin{array}{c}
\scriptstyle{0\le j\le 1}\\
\scriptstyle{i+nj<m+n}\end{array}}b_{ij}x^iy^j,\quad\deg P_1=2n-1-m,\ \deg P_2=
2n-d-m\ .$$ It follows from (\ref{eb20}) that $P_2(0)\ne0$, and hence
$R_1(x,0)$ has the root $x=0$ of some multiplicity $r>0$, while $R_x(x,0)$ has the root $x=0$
with multiplicity at least $r$. However, this contradicts the fact that
$$3nR(x,0)=\left(xR_x(x,y)+nyR_y(x,y)+R_1(x,y)\right)_{y=0}=xR_x(x,0)+R_1(x,0)$$ has the root $x=0$ of
multiplicity $r$. 

The proof of Theorem \ref{tb} is completed.\\


\end{appendix}

{\small
\subsection*{Acknowledgements} The author was supported by the grant no.
176/15 from the Israeli Science Foundation. This work has been completed during the author's stay at the Institute Mittag-Leffler, Stockholm,
and the author is very grateful to IML for hospitality and financial support.}